# $L^1$ BOUNDS IN NORMAL APPROXIMATION

### By Larry Goldstein

#### *University of Southern California*

The zero bias distribution $W^*$ of $W$, defined though the characterizing equation $EWf(W) = \sigma^2 Ef'(W^*)$ for all smooth functions $f$, exists for all $W$ with mean zero and finite variance $\sigma^2$. For $W$ and $W^*$ defined on the same probability space, the $L^1$ distance between $F$, the distribution function of $W$ with $EW = 0$ and $\mathrm{Var}(W) = 1$, and the cumulative standard normal $\Phi$ has the simple upper bound

$$\|F - \Phi\|_1 \leq 2E|W^* - W|.$$

This inequality is used to provide explicit $L^1$ bounds with moderate-sized constants for independent sums, projections of cone measure on the sphere $S(\ell_n^p)$, simple random sampling and combinatorial central limit theorems.

**1. Introduction.** The zero bias transformation and its use in Stein's method [21] for normal approximation was introduced in [10]. There, it was shown that for any mean zero random variable $W$ with finite variance $\sigma^2$, there exists $W^*$ which satisfies

$$(1) \qquad EWf(W) = \sigma^2 Ef'(W^*)$$

for all absolutely continuous $f$ with $E|Wf(W)| < \infty$. We say that such a $W^*$ has the *$W$-zero biased distribution*. Study of the zero bias distribution was motivated by the size bias transformation and Stein's characterization of the normal (see, e.g., [22]), which shows that $Z \sim \mathcal{N}(0, \sigma^2)$ if and only if

$$(2) \qquad EZf(Z) = \sigma^2 Ef'(Z)$$

for all absolutely continuous $f$ with $E|Zf(Z)| < \infty$.

It is helpful to consider the transformation characterized by (1) as a mapping $W \to W^*$ whose domain is the collection of all mean zero distributions

---

with variance $\sigma^2$. From Stein's characterization (2), it is immediate that this transformation has as its unique fixed point the mean zero normal distribution with variance $\sigma^2$. It seems natural, then, that an approximate fixed point of the transformation would be approximately normal and that we can measure the distance of the distribution of $W$ to the normal by the distance between $W$ and its zero bias version $W^*$.

Here, we consider the $L^1$ distance between distribution functions $F$ and $G$ given by

$$\|F - G\|_1 = \int_{-\infty}^{\infty} |F(t) - G(t)|\, dt \tag{3}$$

and known by many names, including Gini's measure of discrepancy, the Kantorovich metric (see [19]), as well as the Wasserstein, Dudley and Fortet–Mourier distance (see, e.g., [3]). If $F$ is the distribution function of a mean zero, variance 1 random variable $W$ and $F^*$ is that of $W^*$ having the $W$-zero biased distribution, Lemma 2.1 of [8] yields that

$$\|F - \Phi\|_1 \leq 2\|F^* - F\|_1, \tag{4}$$

where $\Phi$ is the cumulative distribution function of the standard normal. To bound the right-hand side of (4), it can be convenient to use the dual form of the $L^1$ distance (see [19]) given by

$$\|F - G\|_1 = \inf E|X - Y|, \tag{5}$$

where the infimum is over all couplings of $X$ and $Y$ on a joint space with marginal distributions $F$ and $G$, respectively. Since the dual representation (5) says that $\|F^* - F\|_1$ is upper bounded by $E|W^* - W|$ for any coupling of $W$ and $W^*$, the following result is immediate.

THEOREM 1.1. *Let $W$ be a mean zero, variance 1 random variable with distribution function $F$ and let $W^*$ have the $W$-zero biased distribution and be defined on the same space as $W$. Then, with $\Phi$ the cumulative distribution function of the standard normal,*

$$\|F - \Phi\|_1 \leq 2E|W^* - W|.$$

The goal of this work is to apply Theorem 1.1 to obtain $L^1$ bounds to the normal for a variety of examples and to express the resulting upper bounds as a third-moment-type quantity multiplied by an explicit, moderate constant; in particular, we study sums of independent variables, projections of cone measure, simple random sampling and combinatorial central limit theorems.

In Section 2, we begin by considering the case where $Y = \sum_{i=1}^{n} Y_i$ is the sum of independent mean zero random variables with finite variances $\sigma_i^2 = \mathrm{Var}(Y_i)$, not only to illustrate the method, but also to take advantage of the fact that the particularly simple construction of $Y^*$ in this case allows for the computation of constants in the bound which are explicit functions of the summand distribution. In particular, letting $I$ be an independent random index with distribution

$$(6) \qquad P(I=i) = \frac{\sigma_i^2}{\sum_{j=1}^n \sigma_j^2},$$

the argument proving part (v) of Lemma 2.1 in [10] shows that removing $Y_I$ and replacing it by a variable $Y_I^*$ having the $Y_I$-zero bias distribution, independent of $\{Y_j, j \neq I\}$, gives a variable $Y^*$ with the $Y$-zero bias distribution, that is, that

$$(7) \qquad Y^* = Y - Y_I + Y_I^*$$

has the $Y$-zero biased distribution. We apply this construction and Theorem 1.1 to derive Theorem 2.1 and Corollary 2.1, which yields, for example, that if $F$ is the distribution function of $W = n^{-1/2} \sum_{i=1}^n U_i$, the sum of $n$ i.i.d. variables with the uniform distribution standardized to have mean zero and variance 1, then

$$\|F - \Phi\|_1 \leq \frac{E|U_1|^3}{3\sqrt{n}} \qquad \text{for all } n = 1, 2, \ldots,$$

that is, we obtain a Berry–Esseen type bound, using the $L^1$ metric, with a constant of $1/3$.

In Section 3 we present two constructions of the zero bias distribution $Y^*$ for $Y = \sum_i Y_i$ which can be used in the presence of dependence. Both of these constructions are related to the one used for size biasing which is reviewed in Section 3.1. The first zero bias construction, presented in Section 3.2, can be applied to random vectors $\mathbf{Y} \in \mathbf{R}^n$ which are coordinate symmetric (also called *unconditional*), that is, vectors for which

$$(8) \qquad (Y_1, \ldots, Y_n) =_d (e_1 Y_1, \ldots, e_n Y_n) \qquad \text{for all } (e_1, \ldots, e_n) \in \{-1, 1\}^n.$$

The second construction of $Y^*$, presented in Section 3.3, depends on the existence of an exchangeable pair $(Y', Y'')$ as in Stein [23], whose components have marginal distribution equal to that of $Y$, and which satisfies the linearity condition

$$(9) \qquad E(Y''|Y') = (1-\lambda)Y' \qquad \text{for some } \lambda \in (0, 1).$$

This construction appeared in [10] and was applied in [9] to obtain supremum norm bounds in normal approximation.

The zero bias construction given in Section 3.2 is used in Section 4 to obtain bounds for the normal approximation for one-dimensional projections of the form

$$(10) \qquad Y = \boldsymbol{\theta} \cdot \mathbf{X},$$

where for some $p > 0$, the vector $\mathbf{X} \in \mathbf{R}^n$ has cone measure $\mathcal{C}_p^n$ and $\boldsymbol{\theta} \in \mathbf{R}^n$ is of unit length. To define $\mathcal{C}_p^n$, let

(11)
$$S(\ell_p^n) = \left\{ \mathbf{x} \in \mathbf{R}^n : \sum_{i=1}^n |x_i|^p = 1 \right\} \quad \text{and}$$
$$B(\ell_p^n) = \left\{ \mathbf{x} \in \mathbf{R}^n : \sum_{i=1}^n |x_i|^p \le 1 \right\}.$$

Then, with $\mu^n$ Lebesgue measure in $\mathbf{R}^n$, the cone measure of $A \subset S(\ell_p^n)$ is given by

(12)    $$\mathcal{C}_p^n(A) = \frac{\mu^n([0,1]A)}{\mu^n(B(\ell_p^n))}, \qquad \text{where } [0,1]A = \{ta : a \in A, 0 \le t \le 1\}.$$

Theorem 4.1 provides a normal bound for the projection $Y$ in (10) in terms of explicit and moderate constants and the quantity $\sum_{i=1}^n |\theta_i|^3$ depending on the projection $\boldsymbol{\theta}$. Cone measure, for $p = 1$ and $p = 2$, respectively, includes the special cases of the uniform distribution over the simplex $\sum_{i=1}^n |x_i| = 1$ and the Euclidean sphere $\sum_{i=1}^n x_i^2 = 1$ in $\mathbf{R}^n$. For these two special cases and for $F$ the standardized distribution function of the projection (10), Theorem 4.1 specializes to, respectively,

(13)
$$\|F - \Phi\|_1 \le \frac{9}{\sqrt{2}} \sum_{i=1}^n |\theta_i|^3 + \frac{4}{n+2} \quad \text{and}$$
$$\|F - \Phi\|_1 \le \frac{9}{\sqrt{3}} \sum_{i=1}^n |\theta_i|^3 + \frac{4}{n+2};$$

for $\boldsymbol{\theta} = n^{-1/2}(1, \ldots, 1)$, the sums in (13) are replaced by $n^{-1/2}$.

In Section 5, we turn our attention to simple random sampling of subsets of size $n$ from a set $\mathcal{A}$ of $N$ numerical characteristics, where each subset is selected uniformly, that is, with probability $\binom{N}{n}^{-1}$. The zero bias construction in Section 3.3 is applied in Theorem 5.1 to yield, under some basic nontriviality conditions, the following bound to normality for the distribution function $F$ of the standardized sum of the characteristics in the sample,

$$\|F - \Phi\|_1 \le \frac{4a_3}{\sigma^3} \left( \frac{n(N-n)}{N(N-1)} \right) \left( 1 + \frac{n}{N} \right)^2,$$

where

$$a_3 = \sum_{a \in \mathcal{A}} |a - \bar{a}|^3,$$

$\bar{a}$ is the average of the elements in $\mathcal{A}$ and $\sigma^2$ is the variance of the sum of the sampled characteristics, whose explicit form is given in (70).

In Section 6, we study the accuracy of the normal approximation in the combinatorial central limit theorem. In particular, we apply the zero bias construction in Section 3.3 to variables of the form

$$(14) \qquad Y = \sum_{i=1}^{n} a_{i,\pi(i)},$$

for $n$ a positive integer, $\{a_{i,j}\}_{1 \le i, j \le n}$ the elements of a matrix $A \in \mathbf{R}^{n \times n}$, and $\pi$ a uniformly chosen random permutation on $\mathcal{S}_n$, the symmetric group. Theorem 6.1 yields, for the distribution function $F$ of the standardized variable $Y$ in (14),

$$\|F - \Phi\|_1 \le \frac{a_3}{(n-1)\sigma^3} \left( 16 + \frac{56}{(n-1)} + \frac{8}{(n-1)^2} \right),$$

where

$$a_3 = \sum_{i,j=1}^{n} |a_{ij} - a_{i.} - a_{.j} + a_{..}|^3,$$

$a_{i.}, a_{.j}$ and $a_{..}$ are the averages of $a_{ij}$ over $j, i$ and both $i$ and $j$, respectively, and $\sigma^2$ is the variance of $Y$, whose explicit form is given in (88). When the elements of the population $\mathcal{A}$ or the matrix $A$ behave "typically," the bounds provided by Theorems 5.1 and 6.1 will be of the best order, $n^{-1/2}$.

The zero bias transformation was introduced in [10] to provide smooth function bounds of order $1/n$ for simple random sampling, and the coupling given here in Section 5 for that case is related to the one used there. In [9], the zero bias transformation is used to obtain bounds on the supremum, or $L^\infty$ distance, between the distribution of the sum $Y$ in (14) and the normal, in terms of the maximum of $a_{ij}$; the coupling construction of $W$ to $W^*$ in Section 6 of this paper was first given there. Here, the $L^1$ distance is used and the form of the bounds improved, in that they are expressed in terms of third-moment-type quantities. Also, in [9], supremum norm bounds, again in terms of the maximum of $a_{ij}$, were computed for $Y$ when $\pi$ has a distribution constant on cycle type. The bound (4) was first shown in [8] and applied there to derive the $L^1$ rate of convergence to the normal for hierarchical sequences $X_1, X_2, \ldots$ of random variables whose distributions for some $k \ge 1$ and $f : \mathbf{R}^k \to \mathbf{R}$ satisfy

$$X_{n+1} = f(X_{n,1}, \ldots, X_{n,k}), \qquad n \ge 1,$$

where $X_{n,1}, \ldots, X_{n,k}$ are i.i.d. with distribution equal to that of $X_n$.

**2. Independent variables.** In this section, we demonstrate the application of Theorem 1.1 and the construction (7) to produce $L^1$ bounds with small explicit constants for the distance of the distribution of sums of independent variables to the normal. The utility of Theorem 2.1 below is reflected by the fact that the $L^1$ distance on the left-hand side of (16) requires computation of a convolution, but is bounded on the right by terms which require only the calculation of integrals of the form (3) involving marginal distributions.

The proof of Theorem 2.1 requires the following simple proposition. The first claim is stated in (iii), Section 2.3 of [19]; the second is well known and follows immediately from the dual form (5) of the $L^1$ distance. For $H$ a distribution function on $\mathbf{R}$, let

$$H^{-1}(u) = \sup\{x : H(x) < u\} \qquad \text{for } u \in (0, 1)$$

and let $\mathcal{U}(a, b)$ denote the uniform distribution on $(a, b)$.

PROPOSITION 2.1. *For $F$ and $G$ distribution functions and $U \sim \mathcal{U}(0, 1)$, we have*

$$\|F - G\|_1 = E|F^{-1}(U) - G^{-1}(U)|.$$

*Further, for any $a \geq 0$ and $b \in \mathbf{R}$, where $F_{a,b}$ and $G_{a,b}$ are the distribution functions of $aX + b$ and $aY + b$, respectively,*

$$\|F_{a,b} - G_{a,b}\|_1 = a\|F - G\|_1.$$

Note that one consequence of the proposition is that the $L^1$ distance, as the infimum in (5), can always be achieved. In what follows, we will find it convenient to express relations like the second claim in Proposition 2.1 in a notation where the random variable replaces its distribution function, thus, $\|aX - aY\|_1 = a\|X - Y\|_1$.

THEOREM 2.1. *Let $X_i, i = 1, \ldots, n$, be independent mean zero random variables with variances $\sigma_i^2 = \mathrm{Var}(X_i)$ satisfying $\sum_{i=1}^n \sigma_i^2 = 1$, and let*

$$W = \sum_{i=1}^n X_i.$$

*Then for $F$ the distribution function of $W$ and $\Phi$ that of the standard normal,*

$$\|F - \Phi\|_1 \leq 2E|X_I^* - X_I|, \tag{15}$$

*where $X_i^*$ is any variable having the $X_i$-zero biased distribution, independent of $\{X_j, j \neq i\}$, $i = 1, \ldots, n$, and $I$ is a random index, independent of $\{X_i, X_i^*, i = 1, \ldots, n\}$, with distribution $P(I = i) = \sigma_i^2$.*

*Letting $G_i$ and $G_i^*$ be the distribution functions of $X_i$ and $X_i^*$, respectively, we have*

$$\|F - \Phi\|_1 \le 2 \sum_{i=1}^{n} \sigma_i^2 \|G_i^* - G_i\|_1. \tag{16}$$

*In particular, when $W = n^{-1/2} \sum X_i$ for $X, X_1, \ldots, X_n$ i.i.d. with mean zero, variance 1 and distribution function $G$,*

$$\|F - \Phi\|_1 \le \frac{2}{\sqrt{n}} \|G^* - G\|_1, \tag{17}$$

*and $G^*$, the distribution function of $X^*$, may be given explicitly by*

$$G^*(x) = E[X(X - x)\mathbf{1}(X \le x)]. \tag{18}$$

Proof. The coupling (7) yields $W^* - W = X_I^* - X_I$, with $I$ having distribution as in (6), so (15) follows immediately from Theorem 1.1.

Now, let $U_i, i = 1, \ldots, n$, be a collection of i.i.d. $\mathcal{U}(0, 1)$ variables and set

$$(X_i, X_i^*) = (G_i^{-1}(U_i), (G_i^*)^{-1}(U_i)), \qquad i = 1, \ldots, n;$$

by Proposition 2.1, we have

$$E|X_i^* - X_i| = \|G_i^* - G_i\|_1.$$

Averaging the right-hand side of (15) over $I$ then yields (16) by

$$\|F - \Phi\|_1 \le 2E|X_I^* - X_I| = 2 \sum_{i=1}^{n} \sigma_i^2 E|X_i^* - X_i| = 2 \sum_{i=1}^{n} \sigma_i^2 \|G_i^* - G_i\|_1.$$

When the variables are i.i.d., $\sigma_i^2 = 1/n$, and using Proposition 2.1, the bound becomes

$$2 \sum_{i=1}^{n} \sigma_i^2 \|G_i^* - G_i\|_1 = 2\|G_{1/\sqrt{n}}^* - G_{1/\sqrt{n}}\|_1 = \frac{2}{\sqrt{n}} \|G^* - G\|_1,$$

proving (17).

It is shown in [10] that for $X$ with mean zero and variance 1, the distribution function $G^*$ of $X^*$ is absolutely continuous with respect to Lebesgue measure with density $p^*(x) = -E[X\mathbf{1}(X \le x)]$. Hence, the distribution function of $X^*$ is

$$G^*(x) = -E\left(X \int_{-\infty}^{x} \mathbf{1}(X \le u)\, du\right)$$

$$= -E\left(X \int_{X}^{x} du\, \mathbf{1}(X \le x)\right) = E[X(X - x)\mathbf{1}(X \le x)]. \qquad \square$$

Applying (17) and (18) in particular cases leads to the following corollary.

COROLLARY 2.1.   *Let $B_1, \ldots, B_n$ be i.i.d. Bernoulli variables with success probability $p \in (0,1), q = 1 - p$ and $X_i = (B_i - p)/\sqrt{pq}$. Then for the distribution function $F$ of the sum $W = n^{-1/2} \sum_{i=1}^{n} X_i$, having the standardized binomial $\mathcal{B}(n, p)$ distribution, for every $n = 1, 2, \ldots,$*

$$\|F - \Phi\|_1 \leq \frac{p^2 + q^2}{\sqrt{npq}} = \frac{E|X_1|^3}{\sqrt{n}}$$

[*noting that* $E|X_1|^3 = (p^2 + q^2)/\sqrt{pq}$].

*For $F$ the distribution function of the sum $W = n^{-1/2} \sum_{i=1}^{n} U_i$ of $U_1, \ldots, U_n$, i.i.d. variables with the mean zero, variance 1 uniform distribution $\mathcal{U}[-\sqrt{3}, \sqrt{3}]$, for every $n = 1, 2, \ldots,$*

$$\|F - \Phi\|_1 \leq \frac{\sqrt{3}}{4\sqrt{n}} = \frac{E|X_1|^3}{3\sqrt{n}}$$

(*noting that* $E|X_1|^3 = 3\sqrt{3}/4$).

*If $X$ is any mean zero, variance $\sigma_1^2$ random variable with distribution function $G$ and $Z$ has the $\mathcal{N}(0, \sigma_2^2)$ distribution and is independent of $X$, then when $\sigma_1^2 + \sigma_2^2 = 1$, the distribution function $F$ of the variance 1 sum $W = X + Z$ satisfies*

$$\|F - \Phi\|_1 \leq 2\sigma_1^2 \|G^* - G\|_1.$$

PROOF.   For $X = (B - p)/\sqrt{pq}$, by (18), we have

$$G^*(x) = \frac{pq}{\sqrt{pq}}\left(x + \frac{p}{\sqrt{pq}}\right) \qquad \text{for } x \in \left[\frac{-p}{\sqrt{pq}}, \frac{q}{\sqrt{pq}}\right];$$

that is, $X^*$ is equal in distribution to $(U - p)/\sqrt{pq}$, where $U \sim \mathcal{U}[0, 1]$. Hence, by Proposition 2.1,

$$\|G^* - G\|_1 = \left\|\frac{U - p}{\sqrt{pq}} - \frac{B - p}{\sqrt{pq}}\right\|_1 = \frac{1}{\sqrt{pq}}\|U - B\|_1 = \frac{p^2 + q^2}{2\sqrt{pq}}$$

and the claim now follows by (17) of Theorem 2.1.

For the uniform distribution $\mathcal{U}[-\sqrt{3}, \sqrt{3}]$, (18) yields

$$G^*(x) = -\frac{\sqrt{3}x^3}{36} + \frac{\sqrt{3}x}{4} + \frac{1}{2} \qquad \text{for } x \in [-\sqrt{3}, \sqrt{3}].$$

Now applying (3), we obtain

$$\|G^* - G\|_1 = \frac{\sqrt{3}}{8}.$$

The final claim of the corollary follows from (16) with $n = 2$ and the fact that the normal is a fixed point of the zero bias transformation.  □

Corollary 2.1 yields constants 1 and 1/3 for the standardized Bernoulli, and the Uniform, respectively. Though it is perhaps of greater interest that such constants can be computed explicitly as a function of the underlying distribution, the following proposition gives a bound for the nonidentically distributed case in terms of a universal constant $c_1$, which can be shown to be at most 3. In particular, let

$$(19) \qquad c_1 = \sup \frac{2E|X^* - X|}{E|X|^3},$$

where the supremum is taken over all $X$ with $EX = 0, EX^2 = 1, E|X|^3 < \infty$ and $E|X^* - X| = \|X^* - X\|_1$, that is, with $X^*$ achieving the minimal $L^1$ coupling to $X$.

PROPOSITION 2.2. *For F the distribution function of any variance 1 sum $W = \sum_{i=1}^{n} X_i$ of independent mean zero variables $X_i, i = 1, \ldots, n$,*

$$\|F - \Phi\|_1 \leq c_1 \sum_{i=1}^{n} E|X_i|^3, \qquad \text{where } c_1 \leq 3.$$

PROOF. Let $X$ have mean zero, variance 1, and finite absolute third moment, and let $X^*$ be any variable on the same space as $X$, having the $X$-zero bias distribution. Applying (1), with $f(x) = (1/2)x^2 \operatorname{sgn}(x)$, for which $f'(x) = |x|$, yields

$$E|X^*| = \tfrac{1}{2}E|X|^3.$$

By the triangle inequality and Hölder's inequality, using $EX^2 = 1$ to bound $E|X|$ by $E|X|^3$, we have

$$E|X^* - X| \leq E|X^*| + E|X| \leq \tfrac{1}{2}E|X|^3 + E|X|^3 = \tfrac{3}{2}E|X|^3,$$

yielding $c_1 \leq 3$.

Dropping the requirement that $EX^2 = 1$ in (19), by scaling we have

$$(20) \qquad c_1 = \sup \frac{2\operatorname{Var}(X)E|X^* - X|}{E|X|^3},$$

where the supremum is taken over all $X$ with $EX = 0, 0 < EX^2 < \infty, E|X|^3 < \infty$ and $X^*$ achieving the minimal $L^1$ coupling to $X$.

Now, with $(X_i, X_i^*)$ achieving the minimal $L^1$ coupling for $i = 1, \ldots, n$, (16) and (20) yield

$$\|F - \Phi\|_1 \leq 2\sum_{i=1}^{n} \sigma_i^2 E|X_i^* - X_i|$$

$$= \sum_{i=1}^{n} \left( \frac{2\sigma_i^2 E|X_i^* - X_i|}{E|X_i|^3} \right) E|X_i|^3 \leq c_1 \sum_{i=1}^{n} E|X_i|^3.$$

$\square$

Finally, we remark that as the supremum in (19) is taken over a class of random variables determined by two constraints, the content of [13] and [15] suggests that it may be attained on a three-point distribution.

**3. Coupling constructions.** In this section, we present two constructions which may be used to obtain a variable $Y^*$ having the $Y$-zero bias distribution in the presence of dependence. The first applies when $Y$ is a sum of the components of a coordinate-symmetric vector defined in (8); the second construction uses the exchangeable pair $(Y', Y'')$ of Stein satisfying the linearity condition (9), which first appeared in [10]. We begin by reviewing the construction for size biasing as presented in [11], as both zero bias constructions are related to it.

3.1. *Size biasing.* The zero bias characterization (1) is similar to, and, indeed, was motivated in [10] by, the characterization of the size biased distribution $Y^s$ for a nonnegative random variable $Y$ with finite expectation $\mu$,

$$(21) \qquad\qquad EYf(Y) = \mu Ef(Y^s),$$

holding for all functions $f$ for which $E|Yf(Y)| < \infty$. Under the nontriviality condition $P(Y = 0) < 1$ or, equivalently, the condition $\mu > 0$, the characterization (21) is easily seen to be the same as the more common specification of the size bias distribution $F^s(y)$ as the one which is absolutely continuous with respect to the distribution $F(y)$ of $Y$ with Radon–Nikodym derivative

$$(22) \qquad\qquad \frac{dF^s(y)}{dF(y)} = \frac{y}{\mu}.$$

For the construction of $Y^s$ when $Y = \sum_i Y_i$, the sum of the components of a vector $\mathbf{Y}$ of nonnegative dependent variables with finite means $\mu_i = EY_i$, following [11], we note that for every $i = 1, \ldots, n$, there exists a distribution $\mathbf{Y}^{(i)}$ such that for all functions $f : \mathbf{R}^n \to \mathbf{R}$ for which the expectation on the left-hand side exists,

$$(23) \qquad\qquad EY_if(\mathbf{Y}) = \mu_i Ef(\mathbf{Y}^{(i)});$$

we say that $\mathbf{Y}^{(i)}$ has the $\mathbf{Y}$-size biased distribution in direction $i$. By specializing (23) to the case where $f$ depends only on $Y_i$, we recover (21), showing that $Y_i^{(i)} =_d Y_i^s$, that is, that the $i$th component of $\mathbf{Y}^{(i)}$ has the $Y_i$-size bias distribution.

Without loss of generality, by removing any trivial components of $\mathbf{Y}$ for which $\mu_i = 0$ and lowering the dimension of $\mathbf{Y}$ accordingly, we may express

(23) in the language of (22): denoting the distribution of $\mathbf{Y}$ as $F(\mathbf{y})$, the distribution $F^{(i)}(\mathbf{y})$ of $\mathbf{Y}^{(i)}$ is given by

$$dF^{(i)}(\mathbf{y}) = \frac{y_i}{\mu_i} dF(\mathbf{y}), \tag{24}$$

that is, $\mathbf{Y}^{(i)}$ is absolutely continuous with respect to $\mathbf{Y}$, with Radon–Nikodym derivative $y_i/\mu_i$. Now, as shown in [11], choosing an independent index $I \in \{1, \ldots, n\}$ proportional to the mean of the components of $\mathbf{Y}$, that is, according to the distribution (6), where $\sigma_i^2$ is replaced by $\mu_i$, the variable

$$Y^s = \sum_{j=1}^{n} Y_j^{(I)} \tag{25}$$

has the $Y$-size biased distribution.

Hence, by randomization over $I$, a construction of $\mathbf{Y}^i$ for every $i$ leads to one for $Y^s$. We may accomplish the former as follows. Write the joint distribution of $\mathbf{Y}$ as a product of the marginal distribution of $Y_i$ times the conditional distribution of the remaining variables given $Y_i$,

$$dF(\mathbf{y}) = dF_i(y_i) \, dF(y_1, \ldots, y_{i-1}, y_{i+1} \ldots, y_n | y_i), \tag{26}$$

which gives a factorization of (24) as

$$dF^{(i)}(\mathbf{y}) = dF_i^{(i)}(y_i) \, dF(y_1, \ldots, y_{i-1}, y_{i+1} \ldots, y_n | y_i), \tag{27}$$
$$\text{where } dF_i^{(i)}(y_i) = \frac{y_i}{\mu_i} \, dF_i(y_i).$$

Comparing the relation in (27) between the marginal distributions $F_i^i(y_i)$ and $F_i(y_i)$ with (22) provides an alternate way of seeing that $Y_i^{(i)} =_d Y_i^s$.

The representation (27) says that one may form $\mathbf{Y}^{(i)}$ by first generating $Y_i^{(i)}$ having the $Y_i$-sized biased distribution, and then the remaining variables from their original distribution, conditioned on $y_i$ taking on its newly chosen sized biased value. For $\mathbf{Y}$ already given, a coupling between $Y$ and $Y^s$ can be generated by constructing $Y_i^{(i)}$ and then "adjusting" as necessary the remaining variables $Y_j$, $j \neq i$, so that these have the conditional distribution given $Y_i$ taking on its new value. Typically, the goal is to adjust the variables as little as possible in order to make the resulting bounds to normality small; see [9] and [11] for examples.

In the case where $Y_i$, $i = 1, \ldots, n$, are independent, clearly $Y_j^i =_d Y_j$ for all $j \neq i$. Hence, the construction given above reduces to simply choosing one summand at random with probability proportional to its expectation and replacing it with its biased version. We note that in both zero and size biasing, a finite sum $Y = \sum_i Y_i$ of independent variables is biased by choosing at random and then replacing the randomly chosen variable by a biased version; in size biasing, the variable is chosen proportional to its expectation and in zero biasing, to its variance. The zero bias transformation was so named due to its similarity to size biasing and its application to mean zero random variables.

3.2. *Coordinate symmetric variables.* Of the two zero bias constructions presented here, the one for coordinate symmetric random vectors $\mathbf{Y} \in \mathbf{R}^n$ as defined in (8) is closest to the size biasing construction just described. To begin, note that for all $Y$ such that $EY^2 < \infty$, by replacing the variable $Y$ on the left-hand side of (21) by $Y^2$, we can define the square bias distribution $Y^e$ of $Y$ by the characterization

$$EY^2 f(Y) = EY^2 E f(Y^e)$$

for all functions $f$ for which the expectation of the left-hand side exists. Naturally, when $Y$ has mean zero and variance $\sigma^2$, this identity becomes

$$(28) \qquad EY^2 f(Y) = \sigma^2 E f(Y^e).$$

To make an extension analogous to the one from (21) to (23) for size biasing, let the components of $\mathbf{Y} \in \mathbf{R}^n$ have mean zero and finite variances $\mathrm{Var}(Y_i) = \sigma_i^2$. For such $\mathbf{Y}$, for all $i = 1, \ldots, n$, there exists a distribution $\mathbf{Y}^i$ such that for all functions $f : \mathbf{R}^n \to \mathbf{R}$ for which the expectation of the left-hand side exists,

$$(29) \qquad EY_i^2 f(\mathbf{Y}) = \sigma_i^2 E f(\mathbf{Y}^i);$$

we say that $\mathbf{Y}^i$ has the $\mathbf{Y}$-square biased distribution in direction $i$. By specializing (29) to the case where $f$ depends only on $Y_i$, we recover (28), showing that $Y_i^i =_d Y_i^e$, that is, that the $i$th component of $\mathbf{Y}^i$ has the $Y_i$-square bias distribution.

By removing any component of $\mathbf{Y}$ which is constant and lowering the dimension accordingly, we can assume, without any loss of generality, that each component is nontrivial, that is, that $\sigma_i^2 > 0$ for every $i = 1, \ldots, n$. Parallel to (24) in the case of size biasing, we may now equivalently specify the $\mathbf{Y}^i$ distribution as the one which is absolutely continuous with respect to $\mathbf{Y}$, with

$$(30) \qquad dF^i(\mathbf{y}) = \frac{y_i^2}{\sigma_i^2} dF(\mathbf{y}).$$

Now, let $\mathbf{Y}$ be coordinate-symmetric as defined in (8). Applying (8) marginally and pairwise yields $Y_i =_d -Y_i$ and $(Y_i, Y_j) =_d (-Y_i, Y_j)$ for all $i$ and distinct $i, j$, respectively. Hence, when all components of $\mathbf{Y}$ have finite second moments, taking the following expectation using these distributional equalities yields

$$(31) \qquad EY_i = 0 \qquad \text{for all } i \quad \text{and} \quad EY_i Y_j = 0 \qquad \text{for all } i \neq j.$$

Proposition 3.1 shows how to construct the zero bias distribution $Y^*$ for the sum $Y$ of the components of a coordinate-symmetric vector in terms of $\mathbf{Y}^i$ and a random index in a way that parallels the construction for size biasing given in (25). We let $\mathcal{U}[a, b]$ denote the uniform distribution on $[a, b]$.

PROPOSITION 3.1. *Let* $\mathbf{Y} \in \mathbf{R}^n$ *be a coordinate-symmetric vector as in* (8)*, with* $\operatorname{Var}(Y_i) = \sigma_i^2 \in (0, \infty)$ *for all* $i = 1, 2, \ldots, n$ *and*

$$Y = \sum_{i=1}^{n} Y_i.$$

*Let* $\mathbf{Y}^i, i = 1, \ldots, n,$ *have the "squared bias" distribution given in* (29)*, $I$ be a random index independent of* $\mathbf{Y}$ *and* $\{\mathbf{Y}^i, i = 1, \ldots, n\}$ *with distribution*

$$(32) \qquad P(I = i) = \frac{\sigma_i^2}{\sum_{j=1}^{n} \sigma_j^2}$$

*and* $U \sim \mathcal{U}[-1, 1]$ *be independent of all other variables. Then*

$$(33) \qquad Y^* = U Y_I^I + \sum_{j \neq I} Y_j^I$$

*has the* $Y$-*zero bias distribution.*

PROOF. Let $f$ be an absolutely continuous function with $E|Yf(Y)| < \infty$. Averaging over the index $I$, integrating out the uniform variable $U$ and then applying (29) and (8) to obtain the fourth equality and fifth equalities below, respectively, we have

$$\sigma^2 E f'(Y^*) = \sigma^2 E f'\left(U Y_I^I + \sum_{j \neq I} Y_j^I\right)$$

$$= \sigma^2 \sum_{i=1}^{n} \frac{\sigma_i^2}{\sigma^2} E f'\left(U Y_i^i + \sum_{j \neq i} Y_j^i\right)$$

$$= \sum_{i=1}^{n} \sigma_i^2 E\left(\frac{f(Y_i^i + \sum_{j \neq i} Y_j^i) - f(-Y_i^i + \sum_{j \neq i} Y_j^i)}{2 Y_i^i}\right)$$

$$= \sum_{i=1}^{n} E Y_i\left(\frac{f(Y_i + \sum_{j \neq i} Y_j) - f(-Y_i + \sum_{j \neq i} Y_j)}{2}\right)$$

$$= \sum_{i=1}^{n} E Y_i f\left(Y_i + \sum_{j \neq i} Y_j\right)$$

$$= E Y f(Y).$$

Thus, $Y^*$ has the $Y$-zero bias distribution.   $\square$

The construction for zero biasing implicit in Proposition 3.1 is parallel to the one given in Section 3.1 for size biasing. The factorization (26) suggests that we write (30) as

$$(34) \qquad dF^i(\mathbf{y}) = dF_i^i(y_i)\,dF(y_1,\ldots,y_{i-1},y_{i+1}\ldots,y_n|y_i)$$

$$\text{where } dF_i^i(y_i) = \frac{y_i^2\,dF_i(y_i)}{\sigma_i^2};$$

the relation given in (34) between the marginal distributions $F_i^i(y_i)$ and $F_i(y_i)$ provides an alternate way of seeing that $Y_i^i =_d Y^e$. As for the size biasing construction in Section 3.1, given $\mathbf{Y}$, Proposition 3.1 and (34) now give a coupling between $Y$ and $Y^*$, where an index $I = i$ is chosen with weight proportional to the variance $\sigma_i^2$, the summand $Y_i$ is replaced by $Y_i^i$ having that summand's "square bias" distribution and then multiplied by $U$ and, finally, the remaining variables are adjusted according to their original distribution, given that the $i$th variable takes on the value $Y_i^i$. This construction will be applied in Section 4.

3.3. *Use of the exchangeable pair.* Let $Y$ be a mean zero random variable with finite, nonzero variance. The following description of a coupling of $Y$ to a $Y^*$ having the $Y$-zero biased distribution appears in [10]; its simple proof and some of the consequences below needed for the constructions in Sections 5 and 6 appear in [9].

PROPOSITION 3.2. *Let* $Y',Y''$ *be an exchangeable pair with* $\mathrm{Var}(Y') = \sigma^2 \in (0,\infty)$ *and distribution* $F(y',y'')$ *which satisfies the linearity condition* (9). *Then*

$$(35) \qquad EY' = 0 \quad and \quad E(Y'-Y'')^2 = 2\lambda\sigma^2,$$

*and when* $Y^\dagger, Y^\ddagger$ *have distribution*

$$(36) \qquad dF^\dagger(y',y'') = \frac{(y'-y'')^2}{E(Y'-Y'')^2}\,dF(y',y''),$$

*and* $U \sim \mathcal{U}[0,1]$ *is independent of* $Y^\dagger, Y^\ddagger$, *the variable*

$$Y^* = UY^\dagger + (1-U)Y^\ddagger \qquad \text{has the } Y'\text{-zero biased distribution.}$$

The following construction of $Y^\dagger, Y^\ddagger$ is in the same spirit as the ones given in Sections 3.1 and 3.2. Given $Y'$, first construct $Y''$ close to $Y'$, such that $(Y',Y'')$ is exchangeable and satisfies (9), and use it to form the difference $Y'-Y''$. Then, perhaps independently, construct the parts of $Y^\dagger, Y^\ddagger$ which depend on the "square biased" term $(Y' - Y'')^2$. Finally, construct the remaining parts of $Y^\dagger, Y^\ddagger$ by adjusting the corresponding parts of $Y', Y''$ to have their original joint distribution, given the newly generated variables.

We can describe the constructions used in Sections 5 and 6 in a bit more detail, where the pair $Y', Y''$ is a function of some collection of underlying random variables $\{\xi_\alpha, \alpha \in \mathcal{X}\}$ and an index $\mathbf{I} \subset \mathcal{X}$, possibly random but independent of $\{\xi_\alpha, \alpha \in \mathcal{X}\}$, and the difference $Y' - Y''$ depends only on $\{\xi_\alpha, \alpha \in \mathbf{I}\}$, that is, for some collection of functions $b_\mathbf{I}(\xi_\alpha, \alpha \in \mathbf{i})$,

$$(37) \qquad Y' - Y'' = b_\mathbf{I}(\xi_\alpha, \alpha \in \mathbf{I}).$$

Since one may first generate $\mathbf{I}$, then $\{\xi_\alpha, \alpha \in \mathbf{I}\}$, and finally $\{\xi_\alpha, \alpha \in \mathbf{I}^c\}$ conditional on $\{\xi_\alpha, \alpha \in \mathbf{I}\}$, we may write the joint distribution of all of the variables as

$$(38) \quad dF(\mathbf{i}, \xi_\alpha, \alpha \in \mathcal{X}) = P(\mathbf{I} = \mathbf{i}) \, dF_\mathbf{i}(\xi_\alpha, \alpha \in \mathbf{i}) \, dF_{\mathbf{i}^c|\mathbf{i}}(\xi_\alpha, \alpha \notin \mathbf{i}|\xi_\alpha, \alpha \in \mathbf{i}).$$

Now, consider the distribution $F^\dagger$, which is $F$-square biased by $(y' - y'')^2$:

$$(39) \qquad dF^\dagger(\mathbf{i}, \xi_\alpha, \alpha \in \mathcal{X}) = \frac{(y' - y'')^2}{E(Y' - Y'')^2} \, dF(\mathbf{i}, \xi_\alpha, \alpha \in \mathcal{X}).$$

Using (35) and (37), we obtain

$$2\lambda\sigma^2 = E(Y' - Y'')^2 = Eb_\mathbf{I}^2(\xi_\alpha, \alpha \in \mathbf{I}) = \sum_{\mathbf{i} \subset \mathcal{X}} P(\mathbf{I} = \mathbf{i}) Eb_\mathbf{i}^2(\xi_\alpha, \alpha \in \mathbf{i}),$$

so, in particular, we may define a distribution for an index $\mathbf{I}^\dagger$ with values in subsets of $\mathcal{X}$ by

$$P(\mathbf{I}^\dagger = \mathbf{i}) = \frac{r_\mathbf{i}}{2\lambda\sigma^2} \qquad \text{with } r_\mathbf{i} = P(\mathbf{I} = \mathbf{i}) Eb_\mathbf{i}^2(\xi_\alpha, \alpha \in \mathbf{i}).$$

Hence, substituting (37) and (38) into (39),

$$(40) \quad \begin{aligned} & dF^\dagger(\mathbf{i}, \xi_\alpha, \alpha \in \mathcal{X}) \\ &= \frac{P(\mathbf{I} = \mathbf{i}) b_\mathbf{i}^2(\xi_\alpha, \alpha \in \mathbf{i})}{2\lambda\sigma^2} \, dF_\mathbf{i}(\xi_\alpha, \alpha \in \mathbf{i}) \, dF_{\mathbf{i}^c|\mathbf{i}}(\xi_\alpha, \alpha \notin \mathbf{i}|\xi_\alpha, \alpha \in \mathbf{i}) \\ &= \frac{r_\mathbf{i}}{2\lambda\sigma^2} \frac{b_\mathbf{i}^2(\xi_\alpha, \alpha \in \mathbf{i})}{Eb_\mathbf{i}^2(\xi_\alpha, \alpha \in \mathbf{i})} \, dF_\mathbf{i}(\xi_\alpha, \alpha \in \mathbf{i}) \, dF_{\mathbf{i}^c|\mathbf{i}}(\xi_\alpha, \alpha \notin \mathbf{i}|\xi_\alpha, \alpha \in \mathbf{i}) \\ &= P(\mathbf{I}^\dagger = \mathbf{i}) \, dF_\mathbf{i}^\dagger(\xi_\alpha, \alpha \in \mathbf{i}) \, dF_{\mathbf{i}^c|\mathbf{i}}(\xi_\alpha, \alpha \notin \mathbf{i}|\xi_\alpha, \alpha \in \mathbf{i}), \end{aligned}$$

where

$$dF_\mathbf{i}^\dagger(\xi_\alpha, \alpha \in \mathbf{i}) = \frac{b_\mathbf{i}^2(\xi_\alpha, \alpha \in \mathbf{i})}{Eb_\mathbf{i}^2(\xi_\alpha, \alpha \in \mathbf{i})} \, dF_\mathbf{i}(\xi_\alpha, \alpha \in \mathbf{i}),$$

giving a representation of $dF^\dagger(\mathbf{i}, \xi_\alpha, \alpha \in \mathcal{X})$ parallel to the one for $dF(\mathbf{i}, \xi_\alpha, \alpha \in \mathcal{X})$ in (38). This parallel representation gives a parallel construction as well:

first generate $\mathbf{I}^\dagger$, then $\{\xi_\alpha^\dagger, \alpha \in \mathbf{i}\}$ according to $dF_\mathbf{i}^\dagger$ and finally, $\{\xi_\alpha, \alpha \notin \mathbf{i}\}$ according to $dF_{\mathbf{i}^c|\mathbf{i}}(\xi_\alpha, \alpha \notin \mathbf{i}|\xi_\alpha, \alpha \in \mathbf{i})$.

For the two examples in Sections 5 and 6, the index $\mathbf{I}$ is uniform over some range, so by (40), over that same range, $\mathbf{I}^\dagger$ and $\{\xi_\alpha^\dagger, \alpha \in \mathbf{i}\}$ are jointly drawn from the distribution with proportionality

$$(41) \qquad dF_{\mathbf{i},\xi}^\dagger(\mathbf{i}, \xi_\alpha, \alpha \in i) \sim b_\mathbf{i}^2(\xi_\alpha, \alpha \in \mathbf{i})\, dF_\mathbf{i}(\xi_\alpha, \alpha \in \mathbf{i}).$$

With $\mathbf{I}$ and $\{\xi_\alpha, \alpha \in \mathcal{X}\}$ given, the coupling proceeds by generating $\mathbf{I}^\dagger$ and $\{\xi_\alpha^\dagger, \alpha \in \mathbf{I}^\dagger\}$ according to (41), then adjusting the remaining given variables. For making the bounds small, the goal is to make changes as little as possible, so that the zero biased variable is close to the original.

In Section 5, this procedure results in $S$, a function of the variables which can be kept fixed throughout the construction, and variables $T', T^\dagger$ and $T^\ddagger$ on a joint space such that

$$(42) \qquad Y' = S + T', \qquad Y^\dagger = S + T^\dagger \quad \text{and} \quad Y^\ddagger = S + T^\ddagger,$$

and hence

$$|Y^* - Y'| = |UT^\dagger + (1-U)T^\ddagger - T'|.$$

Here, the underlying variables $\{\xi_\alpha, \alpha \in \mathcal{X}\}$ are $\{X', X'', X_2, \ldots, X_n\}$ and the difference $Y' - Y'' = X' - X''$ so that $\mathbf{I}$ is nonrandom, that is, it indexes the variables $X', X''$ with probability 1, and $b(X', X'') = X' - X''$.

In Section 6, $\{\pi(i), i \in \{1, \ldots, n\}\}$ play the role of $\{\xi_\alpha, \alpha \in \mathcal{X}\}$, $\mathbf{I} = \{I, J\}$ is uniform over all pair of distinct indices in $\{1, \ldots, n\}$ and the difference $Y' - Y''$ is given by

$$(43) \qquad b_{\{i,j\}}(\pi(k), k \in \{i, j\}) = (a_{i,\pi(i)} + a_{j,\pi(j)}) - (a_{i,\pi(j)} + a_{j,\pi(i)}).$$

Note that even when $\mathbf{I}$ is uniformly distributed, the index $\mathbf{I}^\dagger$ need not be; in particular, the distribution (94) given by (41) with $b_\mathbf{i} = b_{\{i,j\}}$ [as in (43)] selects the indices $\mathbf{I}^\dagger = \{I^\dagger, J^\dagger\}$ jointly with their "biased permutation" images $\{K^\dagger, L^\dagger\}$ with probability that preferentially makes the squared difference large. We return to the exchangeable pair construction in Sections 5 and 6.

**4. Projections of cone measure on the sphere $S(\ell_p^n)$.**  In this section, we use the zero biasing construction in Section 3.2 to derive Theorem 4.1, providing bounds to normality for projections $\boldsymbol{\theta} \cdot \mathbf{X}$, where $\mathbf{X} \in \mathbf{R}^n$ has cone measure $\mathcal{C}_p^n$ on the sphere $S(\ell_p^n)$, defined in (12) and (11), respectively, and $\boldsymbol{\theta} \in \mathbf{R}^n$ has unit length. The resulting $L^1$ bound (55) is in terms of explicit small constants [see also (13)] and depends on $\boldsymbol{\theta}$ through the factor $\sum_i |\theta_i|^3$ which yields the best possible rate of $n^{-1/2}$ when the components of $\boldsymbol{\theta}$ are equal.

In the case $p = 2$, cone measure is uniform on the surface of the unit Euclidean sphere in $\mathbf{R}^n$ and [7] shows that the $k$-dimensional projections of $\mathbf{X}$ are close to normal in total variation. The authors of [16] derive normal approximation bounds using Stein's method for random vectors with symmetries in general, including coordinate symmetry, considering the supremum and total variation norm. Studying here the specific instance of cone measure allows for the sharpening of general results to this particular case.

Cone measure is uniform on $S(\ell_p^n)$ only in the cases $p = 1$ and $p = 2$, and the authors of [18] obtain a total variation bound between cone and uniform measure for $p \geq 1$. In some sense, then, the contribution here is related to the central limit problem for convex bodies which strives to quantify when projections of uniform measure on high-dimensional convex bodies have some one-dimensional projection close to normal. A large body of work in this area is generally concerned with the measure of the set of directions on the unit sphere which give rise to approximately normally distributed projections and do not provide bounds in terms of specific projections; see, in particular, [1] and [5] for work continuing that of [24]. In principle, the techniques developed here can be used to shed light on aspects of the central limit theorem for convex bodies; see the remarks at the end of this section.

Let $\mathbf{X} \in \mathbf{R}^n$ be an exchangeable coordinate-symmetric random vector with components having finite second moments and let $\boldsymbol{\theta} \in \mathbf{R}^n$ have unit length. Then, by (31), the projection of $\mathbf{X}$ along the direction $\boldsymbol{\theta}$,

$$Y = \sum_{i=1}^{n} \theta_i X_i,$$

has mean zero and variance $\sigma^2$ equal to the common variance of the components of $\mathbf{X}$. To form $Y^*$ using the construction outlined in Section 3.2, as seen in (34) in particular, requires a vector of random variables to be "adjusted" according to their original distribution, conditional on one coordinate taking on a newly chosen, biased, value. Random vectors which have the "scaling conditional" property in Definition 4.1 can easily be so adjusted. Let $\mathcal{L}(V)$ and $\mathcal{L}(V|X = x)$ denote the distribution of $V$, and the conditional distribution of $V$ given $X = x$, respectively.

DEFINITION 4.1. Let $\mathbf{X} = (X_1, \ldots, X_n)$ be an exchangeable random vector and $\mathcal{D} \subset \mathbf{R}$ the support of the distribution of $X_1$. If there exists a function $g : \mathcal{D} \to \mathbf{R}$ such that $P(g(X_1) = 0) = 0$ and

$$(44) \quad \mathcal{L}(X_2, \ldots, X_n | X_1 = a) = \mathcal{L}\left(\frac{g(a)}{g(X_1)}(X_2, \ldots, X_n)\right) \qquad \text{for all } a \in \mathcal{D},$$

then we say that $\mathbf{X}$ is scaling $g$-conditional or, more simply, scaling-conditional.

Proposition 4.1 is an application of Theorem 1.1 to projections of scaling-conditional vectors.

PROPOSITION 4.1.  *Let $\mathbf{X} \in \mathbf{R}^n$ be an exchangeable, coordinate-symmetric and scaling $g$-conditional random vector with finite second moments and, with $\boldsymbol{\theta} \in \mathbf{R}^n$ of unit length, set*

$$Y = \sum_{i=1}^{n} \theta_i X_i, \qquad \sigma^2 = \mathrm{Var}(Y) \quad and \quad F(x) = P(Y/\sigma \leq x).$$

*Then any construction of $(\mathbf{X}, X_i^i)$ on a joint space for each $i = 1, \ldots, n$ with $X_i^i$ having the $X_i$-square biased distribution provides the upper bound*

$$(45) \qquad \|F - \Phi\|_1 \leq \frac{2}{\sigma} E \left| \theta_I (U X_I^I - X_I) + \left( \frac{g(X_I^I)}{g(X_I)} - 1 \right) \sum_{j \neq I} \theta_j X_j \right|,$$

*where $P(I = i) = \theta_i^2$, $U \sim \mathcal{U}[-1, 1]$ and $I$ and $U$ are independent of each other and of the remaining variables.*

PROOF.  For all $i = 1, \ldots, n$, since $\mathbf{X}$ is scaling $g$-conditional, given $\mathbf{X}$ and $X_i^i$, the vector

$$\mathbf{X}^i = \left( \frac{g(X_i^i)}{g(X_i)} X_1, \ldots, \frac{g(X_i^i)}{g(X_i)} X_{i-1}, X_i^i, \frac{g(X_i^i)}{g(X_i)} X_{i+1}, \ldots, \frac{g(X_i^i)}{g(X_i)} X_n \right)$$

has the $\mathbf{X}$-square bias distribution in direction $i$ as given in (29); in particular, for every $h$ for which the expectation on the left-hand side below exists,

$$(46) \qquad E X_i^2 h(\mathbf{X}) = E X_i^2 E h(\mathbf{X}^i).$$

We now apply Proposition 3.1 to $\mathbf{Y} = (\theta_1 X_1, \ldots, \theta_n X_n)$. First, the coordinate symmetry of $\mathbf{Y}$ follows from that of $\mathbf{X}$. Next, we claim

$$\mathbf{Y}^i = (\theta_1 X_1^i, \ldots, \theta_n X_n^i)$$

has the $\mathbf{Y}$-square bias distribution in direction $i$. Given $f$, applying (46) with

$$h(\mathbf{X}) = f(\theta_1 X_1, \ldots, \theta_n X_n) = f(\mathbf{Y})$$

and then multiplying both sides by $\theta_i^2$ yields

$$E \theta_i^2 X_i^2 f(\mathbf{Y}) = E \theta_i^2 X_i^2 E f(\mathbf{Y}^i) \quad \text{or} \quad E Y_i^2 f(\mathbf{Y}) = E Y_i^2 E f(\mathbf{Y}^i).$$

Finally, since $\mathbf{X}$ is exchangeable, the variance of $Y_i$ is proportional to $\theta_i^2$ and the distribution of $I$ in (32) specializes to the one claimed.

Now, (33) of Proposition 3.1 yields, with $Y^*$ having the $Y$-zero biased distribution,

$$Y^* - Y = UY_I^I + \sum_{j \neq I} Y_j^I - \sum_{i=1}^n Y_i$$

$$= U\theta_I X_I^I + \sum_{j \neq I} \theta_j X_j^I - \sum_{j=1}^n \theta_j X_j$$

$$= \theta_I(UX_I^I - X_I) + \sum_{j \neq I} \theta_j(X_j^I - X_j)$$

$$= \theta_I(UX_I^I - X_I) + \sum_{j \neq I} \theta_j\left(\frac{g(X_I^I)}{g(X_I)} - 1\right)X_j$$

$$= \theta_I(UX_I^I - X_I) + \left(\frac{g(X_I^I)}{g(X_I)} - 1\right)\sum_{j \neq I} \theta_j X_j.$$

The proof is completed by dividing both sides by $\sigma$, noting that $Y^*/\sigma = (Y/\sigma)^*$, and invoking Theorem 1.1. $\quad\square$

Proposition 4.2 shows that Proposition 4.1 can be applied when $\mathbf{X}$ has cone measure. We denote the Gamma and Beta distributions with parameters $\alpha, \beta$ as $\Gamma(\alpha, \beta)$ and $B(\alpha, \beta)$, respectively, and the Gamma function at $x$ by $\Gamma(x)$.

PROPOSITION 4.2. *Let $\mathcal{C}_p^n$ denote cone measure as given in (12) for some $p > 0$.*

1. *Cone measure $\mathcal{C}_p^n$ is exchangeable and coordinate-symmetric. For $\{G_j, \varepsilon_j, j = 1, \ldots, n\}$ independent variables with $G_j \sim \Gamma(1/p, 1)$ and $\varepsilon_j$ taking values $-1$ and $+1$ with equal probability,*

$$(47) \quad \mathbf{X} = \left(\varepsilon_1\left(\frac{G_1}{G_{1,n}}\right)^{1/p}, \ldots, \varepsilon_n\left(\frac{G_n}{G_{1,n}}\right)^{1/p}\right) \sim \mathcal{C}_p^n, \qquad \text{where } G_{a,b} = \sum_{i=a}^b G_i.$$

2. *The common marginal distribution $X_i$ of cone measure is characterized by*

$$X_i =_d -X_i \quad \text{and} \quad |X_i|^p \sim B(1/p, (n-1)/p),$$

   *and the variance $\sigma_{n,p}^2 = \text{Var}(X_i)$ is given by*

$$(48) \quad \sigma_{n,p}^2 = \frac{\Gamma(3/p)\Gamma(n/p)}{\Gamma(1/p)\Gamma((n+2)/p)} \qquad \text{satisfying } \lim_{n \to \infty} n^{2/p}\sigma_{n,p}^2 = \frac{\Gamma(3/p)}{\Gamma(1/p)}.$$

3. *The square bias distribution $X_i^i$ of $X_i$ is characterized by*

$$X_i^i =_d -X_i^i \quad and \quad |X_i^i|^p \sim B(3/p, (n-1)/p). \tag{49}$$

*In particular, letting $\{G_j, G_j', \varepsilon_j, j = 1, \ldots, n\}$ be independent variables with $G_j \sim \Gamma(1/p, 1)$, $G_j' \sim \Gamma(2/p, 1)$ and $\varepsilon_j$ taking values $-1$ and $+1$ with equal probability, for each $i = 1, \ldots, n$, a construction of $(\mathbf{X}, X_i^i)$ on a joint space is given by the representation of $\mathbf{X}$ in (47) along with*

$$X_i^i = \varepsilon_i \left( \frac{G_i + G_i'}{G_{1,n} + G_i'} \right)^{1/p}. \tag{50}$$

*The mean $m_{n,p} = E|X_i^i|$ for all $i = 1, \ldots, n$ is given by*

$$m_{n,p} = \frac{\Gamma(4/p)\Gamma((n+2)/p)}{\Gamma(3/p)\Gamma((n+3)/p)} \tag{51}$$

*and satisfies*

$$\lim_{n \to \infty} n^{1/p} m_{n,p} = \frac{\Gamma(4/p)}{\Gamma(3/p)} \quad and \quad m_{n,p} \le \left( \frac{3}{n+2} \right)^{1/(p \vee 1)}. \tag{52}$$

4. *Cone measure $\mathcal{C}_p^n$ is scaling $(1 - |x|^p)^{1/p}$-conditional.*

The proof of Proposition 4.2 is deferred to the end of this section. Before proceeding to Theorem 4.1, we remind the reader of the following known facts about the Gamma and Beta distributions; see [4], Theorem 1.2.3 for the case $n = 2$ of the first claim, the extension to general $n$ and the following claim being straightforward. For $\gamma_i \sim \Gamma(\alpha_i, 1)$, $i = 1, \ldots, n$, independent and $\alpha_i > 0$,

$$\gamma_1 + \gamma_2 \sim \Gamma(\alpha_1 + \alpha_2, 1), \qquad \frac{\gamma_1}{\gamma_1 + \gamma_2} \sim B(\alpha_1, \alpha_2), \quad and \tag{53}$$

$$\left( \frac{\gamma_1}{\sum_{i=1}^n \gamma_i}, \ldots, \frac{\gamma_n}{\sum_{i=1}^n \gamma_i} \right) \quad and \quad \sum_{i=1}^n \gamma_i \quad are\ independent;$$

the Beta distribution $B(\alpha, \beta)$ has density

$$p_{\alpha,\beta}(u) = \frac{\Gamma(\alpha + \beta)}{\Gamma(\alpha)\Gamma(\beta)} u^{\alpha-1} (1-u)^{\beta-1} \mathbf{1}_{u \in [0,1]}$$

$$and\ \kappa > 0\ moments \quad \frac{\Gamma(\alpha + \kappa)\Gamma(\alpha + \beta)}{\Gamma(\alpha + \beta + \kappa)\Gamma(\alpha)}. \tag{54}$$

THEOREM 4.1. *Let* $\mathbf{X}$ *have cone measure* $\mathcal{C}_p^n$ *on the sphere* $S(\ell_p^n)$ *for some* $p > 0$ *and let*

$$Y = \sum_{i=1}^n \theta_i X_i$$

*be the one-dimensional projection of* $\mathbf{X}$ *along the direction* $\boldsymbol{\theta} \in \mathbf{R}^n$ *with* $\|\boldsymbol{\theta}\| = 1$. *Then with* $\sigma_{n,p}^2 = \mathrm{Var}(X_1)$ *and* $m_{n,p} = E|X_1^1|$ *given in* (48) *and* (51), *respectively, and* $F$ *the distribution function of the normalized sum* $W = Y/\sigma_{n,p}$,

$$(55) \qquad \|F - \Phi\|_1 \le 3\left(\frac{m_{n,p}}{\sigma_{n,p}}\right)\sum_{i=1}^n |\theta_i|^3 + \left(\frac{1}{p} \vee 1\right)\frac{4}{n+2},$$

*where* $\Phi$ *is the cumulative distribution function of the standard normal.*

We note that by the limits in (48) and (52), the constant $m_{n,p}/\sigma_{n,p}$ that multiplies the sum in the bound (55) is of the order of a constant, with asymptotic value

$$\lim_{n \to \infty} \frac{m_{n,p}}{\sigma_{n,p}} = \frac{\Gamma(4/p)\sqrt{\Gamma(1/p)}}{\Gamma(3/p)^{3/2}}.$$

Since, for $\boldsymbol{\theta} \in \mathbf{R}^n$ with $\|\boldsymbol{\theta}\| = 1$, we have

$$\sum |\theta_i|^3 \ge \frac{1}{\sqrt{n}},$$

the second term in (55) is always of smaller order than the first, so the decay rate of the bound to zero is determined by $\sum_i |\theta_i|^3$. The minimal rate $1/\sqrt{n}$ is achieved when $\theta_i = 1/\sqrt{n}$.

In the special cases $p = 1$ and $p = 2$, $\mathcal{C}_p^n$ is uniform on the simplex $\sum_{i=1}^n |x_i| = 1$ and the unit Euclidean sphere $\sum_{i=1}^n x_i^2 = 1$, respectively. By (48) and (51) for $p = 1$,

$$\sigma_{n,1}^2 = \frac{2}{n(n+1)} \quad \text{and} \quad m_{n,1} = \frac{3}{n+2},$$

and, using (52) for $p = 2$,

$$\sigma_{n,2}^2 = \frac{1}{n} \quad \text{and} \quad m_{n,2} \le \sqrt{\frac{3}{n+2}};$$

these relations yield

$$\frac{m_{n,1}}{\sigma_{n,1}} = 3\sqrt{\frac{n(n+1)}{2(n+2)^2}} \le \frac{3}{\sqrt{2}} \quad \text{and} \quad \frac{m_{n,2}}{\sigma_{n,2}} \le \sqrt{\frac{3n}{n+2}} \le \sqrt{3}.$$

Substituting into (55) now gives the claim (13).

PROOF OF THEOREM 4.1. Using Proposition 4.2, we apply Proposition 4.1 for $\mathbf{X}$ with $g(x) = (1 - |x|^p)^{1/p}$ and the joint construction of $(\mathbf{X}, X_i^i)$ given in item 3.

Using the triangle inequality on (45) yields the upper bound

$$(56) \qquad \frac{2}{\sigma_{n,p}} \left( E|\theta_I (UX_I^I - X_I)| + E \left| \left( \frac{g(X_I^I)}{g(X_I)} - 1 \right) \sum_{j \neq I} \theta_j X_j \right| \right).$$

For $X$ with the common marginal of $\mathbf{X}$, we have

$$E \left| \frac{X}{\sigma_{n,p}} \right| \leq \left( E \left| \frac{X}{\sigma_{n,p}} \right|^2 \right)^{1/2} = 1 \leq \left( E \left| \frac{X}{\sigma_{n,p}} \right|^3 \right)^{1/3} \leq \left( E \left| \frac{X}{\sigma_{n,p}} \right|^3 \right)$$

which, with $X^1$ having the square bias distribution of $X$, implies that

$$E|X| \leq \frac{E|X|^3}{\sigma_{n,p}^2} = E|X^1|.$$

Bounding the first term in (56) by applying the triangle inequality, using the fact that $U$ is independent of $I$ and $X_I^I$, $E|U| = 1/2$ and $P(I = i) = \theta_i^2$ yields

$$(57) \qquad \begin{aligned} E|\theta_I|(|UX_I^I| + |X_I|) &= E|\theta_I|(\tfrac{1}{2}|X_I^I| + |X_I|) = E \sum_{i=1}^{n} |\theta_i|^3 (\tfrac{1}{2}|X_i^i| + |X_i|) \\ &\leq \tfrac{3}{2} \sum_{i=1}^{n} |\theta_i|^3 E|X_i^i| = \tfrac{3}{2} m_{n,p} \sum_{i=1}^{n} |\theta_i|^3. \end{aligned}$$

Now, averaging the second term in (56) over the distribution of $I$ yields

$$(58) \qquad E \left| \left( \frac{g(X_I^I)}{g(X_I)} - 1 \right) \sum_{j \neq I} \theta_j X_j \right| = \sum_{i=1}^{n} E \left| \left( \frac{g(X_i^i)}{g(X_i)} - 1 \right) \sum_{j \neq i} \theta_j X_j \right| \theta_i^2.$$

Using (47), (50) and $g(x) = (1 - |x|^p)^{1/p}$, we have

$$(59) \qquad \frac{g(X_i^i)}{g(X_i)} - 1 = \left( \frac{G_{1,n}}{G_{1,n} + G_i'} \right)^{1/p} - 1.$$

The variable $G_i'$ and, by (53), the sum $G_{1,n}$ are independent of $X_1, \ldots, X_n$; hence, the term (59) is independent of the sum it multiplies in (58) and therefore equals

$$(60) \qquad \sum_{i=1}^{n} E \left| \frac{g(X_i^i)}{g(X_i)} - 1 \right| E \left| \sum_{j \neq i} \theta_j X_j \right| \theta_i^2.$$

To bound the first expectation in (60), since $G_{1,n}/(G_{1,n} + G_i') \sim B(n/p, 2/p)$, we have

$$(61) \qquad E\left|\frac{g(X_i^i)}{g(X_i)} - 1\right| = E\left(1 - \left(\frac{G_{1,n}}{G_{1,n} + G_i'}\right)^{1/p}\right) \leq \left(\frac{1}{p} \vee 1\right)\frac{2}{n+2}$$

since for $p \geq 1$, using (54) with $\kappa = 1$,

$$E\left(1 - \left(\frac{G_{1,n}}{G_{1,n} + G_i'}\right)^{1/p}\right) \leq E\left(1 - \left(\frac{G_{1,n}}{G_{1,n} + G_i'}\right)\right) = 1 - \frac{n/p}{(n+2)/p} = \frac{2}{n+2},$$

while for $0 < p < 1$, using Jensen's inequality and the fact that $(1-x)^{1/p} \geq 1 - x/p$ for $x \leq 1$,

$$E\left(1 - \left(\frac{G_{1,n}}{G_{1,n} + G_i'}\right)^{1/p}\right) \leq 1 - \left(E\frac{G_{1,n}}{G_{1,n} + G_i'}\right)^{1/p}$$

$$= 1 - \left(\frac{n}{n+2}\right)^{1/p} \leq \frac{2}{p(n+2)}.$$

We may bound the second expectation in (60) by $\sigma_{n,p}$ since

$$\left(E\left|\sum_{j \neq i} \theta_j X_j\right|\right)^2 \leq E\left(\sum_{j \neq i} \theta_j X_j\right)^2 = \text{Var}\left(\sum_{j \neq i} \theta_j X_j\right) = \sigma_{n,p}^2 \sum_{j \neq i} \theta_j^2 \leq \sigma_{n,p}^2.$$

Neither this bound nor the bound (61) depends on $i$, so substituting them into (60) and summing over $i$, again using $\sum_i \theta_i^2 = 1$, yields

$$(62) \qquad \sum_{i=1}^n E\left|\frac{g(X_i^i)}{g(X_i)} - 1\right| E\left|\sum_{j \neq i} \theta_j X_j\right| \theta_i^2 \leq \sigma_{n,p}\left(\frac{1}{p} \vee 1\right)\frac{2}{n+2}.$$

Adding (57) and (62) and multiplying by $2/\sigma_{n,p}$ in accordance with (45) yields (55). $\square$

PROOF OF PROPOSITION 4.2.

1. For $A \subset S(\ell_p^n)$, $\mathbf{e} = (e_1, \ldots, e_n) \in \{-1, 1\}^n$ and a permutation $\pi \in \mathcal{S}_n$, let

$$A_{\mathbf{e}} = \{\mathbf{x} : (e_1 x_1, \ldots, e_n x_n) \in A\} \quad \text{and} \quad A_\pi = \{\mathbf{x} : (x_{\pi(1)}, \ldots, x_{\pi(n)}) \in A\}.$$

By the properties of Lebesgue measure, $\mu^n([0,1]A_{\mathbf{e}}) = \mu^n([0,1]A_\pi) = \mu^n([0,1]A)$, so by (12), cone measure is coordinate symmetric and exchangeable.

Next, [20], for instance, shows that

$$(63) \qquad (|X_1|, \ldots, |X_n|) =_d \left(\left(\frac{G_1}{G_{1,n}}\right)^{1/p}, \ldots, \left(\frac{G_n}{G_{1,n}}\right)^{1/p}\right).$$

Letting $C$ and $C_{\mathbf{e}}$ be the distribution functions of $\mathbf{X} \sim \mathcal{C}_p^n$ and $(e_1 X_1, \ldots, e_n X_n)$, respectively, the coordinate symmetry of $\mathbf{X}$ implies that

$$C(\mathbf{x}) = C_{\mathbf{e}}(\mathbf{x}) \qquad \text{for all } \mathbf{e} \in \{-1, 1\}^n,$$

so averaging yields

$$C(\mathbf{x}) = \frac{1}{2^n} \sum_{\mathbf{e} \in \{-1, 1\}^n} C_{\mathbf{e}}(\mathbf{x}).$$

Therefore, for $\varepsilon_i$, $i = 1, \ldots, n$, i.i.d. variables taking the values 1 and $-1$ with probability $1/2$, we conclude that $\mathbf{X} =_d (\varepsilon_1 X_1, \ldots, \varepsilon_n X_n) =_d (\varepsilon_1 |X_1|, \ldots, \varepsilon_n |X_n|)$. Combining this fact with (63) yields (47).

2. Applying the coordinate symmetry of $\mathbf{X}$ coordinatewise gives $X_i =_d -X_i$ and (63) yields $|X_i|^p = G_i / G_{1,n}$, which has the claimed Beta distribution, by (53). As $EX_i = 0$, we have

$$(64) \qquad \mathrm{Var}(X_i) = EX_i^2 = E(|X_i|^p)^{2/p}$$

and the variance claim in (48) follows from (54) for $\alpha = 1/p, \beta = (n-1)/p$ and $\kappa = 2/p$. The limit in (48) follows from the fact that for all $n, x > 0$,

$$(65) \qquad \lim_{n \to \infty} \frac{n^x \Gamma(n)}{\Gamma(n + x)} = 1,$$

which can be shown using Stirling's formula.

3. If $X$ is symmetric with variance $\sigma^2$ and $X^1$ has the $X$-square bias density, then for all odd functions $f$, since $-X^2 f(X) =_d X^2 f(X)$,

$$Ef(-X^1) = \frac{EX^2 f(-X)}{\sigma^2} = \frac{E(-X^2 f(X))}{\sigma^2} = \frac{E(X^2 f(X))}{\sigma^2} = Ef(X^1),$$

showing that $X^1$ is symmetric.

From (54) and a change of variables, $X$ satisfies $|X|^p \sim B(\alpha/p, \beta/p)$ if and only if the density $p_{|X|}(u)$ of $|X|$ is

$$(66) \qquad p_{|X|}(u) = \frac{p \Gamma((\alpha + \beta)/p)}{\Gamma(\alpha/p) \Gamma(\beta/p)} u^{\alpha-1} (1 - u^p)^{\beta/p-1} \mathbf{1}_{u \in [0,1]}.$$

Hence, since $|X_i|^p \sim B(1/p, (n-1)/p)$ by item 2, the density $p_{|X_i|}(u)$ of $|X_i|$ is

$$p_{|X_i|}(u) = \frac{p \Gamma(n/p)}{\Gamma(1/p) \Gamma((n-1)/p)} (1 - u^p)^{(n-1)/p-1} \mathbf{1}_{u \in [0,1]}.$$

Multiplying by $u^2$ and renormalizing produces the $|X_i^i|$ density

$$(67) \qquad \begin{aligned} p_{|X_i^i|}(u) &= \frac{u^2 p_{|X|}(u)}{EX_i^2} \\ &= \frac{p \Gamma((n+2)/p)}{\Gamma(3/p) \Gamma((n-1)/p)} u^2 (1 - u^p)^{(n-1)/p-1} \mathbf{1}_{u \in [0,1]}, \end{aligned}$$

and comparing (67) to (66) shows the second claim in (49). The representation (50) now follows from (53) and the symmetry of $X_i^i$.

As in (64), the moment formula (51) follows from (54) for $\alpha = 3/p, \beta = (n-1)/p$ and $\kappa = 1/p$, and the limit in (52) follows by (65). Regarding the last claim in (52), for $p \geq 1$, Hölder's inequality gives

$$m_{n,p} = E|X^1| \leq (E|X^1|^p)^{1/p} = \left(\frac{3}{n+2}\right)^{1/p},$$

while for $0 < p < 1$, we have

$$m_{n,p} = E|X^1| = E\left(\frac{G_i + G_i'}{G_{1,n} + G_i'}\right)^{1/p} \leq E\left(\frac{G_i + G_i'}{G_{1,n} + G_i'}\right) = \frac{3}{n+2}.$$

4. We consider the conditional distribution on the left-hand side of (44) and use the representation (and notation $G_{a,b}$) given in (47). The second equality below follows from the coordinate symmetry of $\mathbf{X}$, and the fourth follows since we may replace $G_{1,n}$ by $G_{2,n}/(1 - |a|^p)$ on the conditioning event. Further, using the notation $a\mathcal{L}(V)$ for the distribution of $aV$, we have

$$
\begin{aligned}
&\mathcal{L}(X_2, \ldots, X_n | X_1 = a) \\
&\quad = \mathcal{L}\left(\varepsilon_2 \left(\frac{G_2}{G_{1,n}}\right)^{1/p}, \ldots, \varepsilon_n \left(\frac{G_n}{G_{1,n}}\right)^{1/p} \middle| \varepsilon_1 \left(\frac{G_1}{G_{1,n}}\right)^{1/p} = a\right) \\
&\quad = \mathcal{L}\left(\varepsilon_2 \left(\frac{G_2}{G_{1,n}}\right)^{1/p}, \ldots, \varepsilon_n \left(\frac{G_n}{G_{1,n}}\right)^{1/p} \middle| \left(\frac{G_1}{G_{1,n}}\right)^{1/p} = |a|\right) \\
&\quad = \mathcal{L}\left(\varepsilon_2 \left(\frac{G_2}{G_{1,n}}\right)^{1/p}, \ldots, \varepsilon_n \left(\frac{G_n}{G_{1,n}}\right)^{1/p} \middle| \frac{G_{2,n}}{G_{1,n}} = 1 - |a|^p\right) \\
&\quad = (1 - |a|^p)^{1/p} \mathcal{L}\left(\varepsilon_2 \left(\frac{G_2}{G_{2,n}}\right)^{1/p}, \ldots, \varepsilon_n \left(\frac{G_n}{G_{2,n}}\right)^{1/p} \middle| \frac{G_{2,n}}{G_{1,n}} = 1 - |a|^p\right) \\
&\quad = (1 - |a|^p)^{1/p} \mathcal{L}\left(\varepsilon_2 \left(\frac{G_2}{G_{2,n}}\right)^{1/p}, \ldots, \varepsilon_n \left(\frac{G_n}{G_{2,n}}\right)^{1/p} \middle| \frac{G_1}{G_{1,n}} = |a|^p\right) \\
&\quad = (1 - |a|^p)^{1/p} \mathcal{L}\left(\varepsilon_2 \left(\frac{G_2}{G_{2,n}}\right)^{1/p}, \ldots, \varepsilon_n \left(\frac{G_n}{G_{2,n}}\right)^{1/p}\right) \\
&\quad = g(a)\mathcal{C}_p^{n-1}.
\end{aligned}
$$
(68)

In the penultimate step, we remove the conditioning on $G_1/G_{1,n}$ since (53) and the independence of $G_1$ from all other variables gives that

$$\left(\frac{G_2}{G_{2,n}}, \ldots, \frac{G_n}{G_{2,n}}\right) \qquad \text{is independent of } (G_1, G_{2,n})$$

and so, in particular, is independent of $G_1/(G_1 + G_{2,n}) = G_1/G_{1,n}$.

Regarding the right-hand side of (44), using $1 - |X_1|^p = \sum_{i=2}^{n} |X_i|^p$ and the representation (47), we obtain

$$
\begin{aligned}
g(a)(X_2, \ldots, X_n)/g(X_1) &= g(a)\left(\frac{(X_2, \ldots, X_n)}{(|X_2|^p + \cdots + |X_n|^p)^{1/p}}\right) \\
&= g(a)\left(\frac{(\varepsilon_2(G_2/G_{1,n})^{1/p}, \ldots, \varepsilon_n(G_n/G_{1,n})^{1/p})}{((G_2/G_{1,n}) + \cdots + (G_n/G_{1,n}))^{1/p}}\right) \\
&= g(a)\left(\frac{(\varepsilon_2 G_2^{1/p}, \ldots, \varepsilon_n G_n^{1/p})}{(G_2 + \cdots + G_n)^{1/p}}\right) \\
&= g(a)\left(\varepsilon_2\left(\frac{G_2}{G_{2,n}}\right)^{1/p}, \ldots, \varepsilon_n\left(\frac{G_n}{G_{2,n}}\right)^{1/p}\right) \\
&=_d g(a)\mathcal{C}_p^{n-1},
\end{aligned}
$$

matching the distribution (68). □

In principle, Proposition 3.1 can be applied in conjunction with Theorem 1.1 for any coordinate-symmetric vector where one can construct a coupling between the marginal variables and their square biased versions, and where conditional distributions such as the one on the left-hand side of (44) can be handled. For **X** having the uniform distribution over a convex body symmetric to the coordinate planes, the conditional distributions of interest are uniform over the intersection of the body with the hyperplanes $X_i = a$. The marginal coupling appears to be more elusive, but may be especially tractable when the body has some particular shapes.

**5. Simple random sampling.** We provide an $L^1$ bound for the error in the normal approximation of the sum

$$
Y = \sum_{i=1}^{n} X_i \tag{69}
$$

of a simple random sample of size $n$ from a set $\mathcal{A}$ of $N$ real numbers, not all equal. It is straightforward to verify that $Y$ has mean $\mu$ and variance $\sigma^2$ given by

$$
\mu = n\bar{a} \quad \text{and} \quad \sigma^2 = \frac{n(N-n)}{N(N-1)} \sum_{a \in \mathcal{A}} (a - \bar{a})^2 \qquad \text{where} \quad \bar{a} = \frac{1}{N} \sum_{a \in \mathcal{A}} a. \tag{70}
$$

The bound below depends also on $a_3$, the third-moment-type quantity given by

$$
a_3 = \sum_{a \in \mathcal{A}} |a - \bar{a}|^3.
$$

THEOREM 5.1. *Let $\{X_1, \ldots, X_n\}$ be a simple random sample of size $n$ from a set $\mathcal{A}$ of $N$ real numbers, not all equal, with $n$ and $N$ satisfying*

$$2 < n < N - 1. \tag{71}$$

*Then, with the sum $Y$ given by* (69), *the distribution function $F$ of the standardized variable $W = (Y - \mu)/\sigma$ satisfies*

$$\|F - \Phi\|_1 \leq \frac{4a_3}{\sigma^3} \left( \frac{n(N-n)}{N(N-1)} \right) \left( 1 + \frac{n}{N} \right)^2.$$

Using $n/N \leq 1$, we see that the theorem provides the "universal" upper bound $16a_3/\sigma^3$, although if the sampling fraction $n/N$ is close to $1/2$, the bound improves substantially, close to $2.25a_3/\sigma^3$.

Since $W$ and $a_3/\sigma^3$ are invariant upon replacing $a$ by $(a - \bar{a})/\sqrt{\sum_{b \in \mathcal{A}} (b - \bar{a})^2}$, we may assume below, without loss of generality, that the collection $\mathcal{A}$ satisfies

$$\sum_{a \in \mathcal{A}} a = 0 \quad \text{and} \quad \sum_{a \in \mathcal{A}} a^2 = 1. \tag{72}$$

If we consider a sequence $\mathcal{A}_N$ of collections of $N$ numbers, not all equal, then the bound will be of (the best) order $1/\sqrt{N}$ as $N \to \infty$ if the deviations $a - \bar{a}$, $a \in \mathcal{A}_N$ are comparable and the sampling fraction $n/N$ is bounded away from 0 and 1; in particular, under (72), $\sigma^2$ will be of order 1, the deviations of order $1/\sqrt{N}$ and $a_3$ (and therefore the bound) of order $1/\sqrt{N}$.

PROOF. By (72),

$$\sigma^2 = \frac{n(N-n)}{N(N-1)} \quad \text{and} \quad a_3 = \sum_{a \in \mathcal{A}} |a|^3,$$

so it suffices to prove that

$$\|F - \Phi\|_1 \leq \frac{4a_3}{\sigma} \left( 1 + \frac{n}{N} \right)^2. \tag{73}$$

Since distinct labels may be appended to the elements of $\mathcal{A}$, say as a second coordinate which is neglected when taking sums, we may assume that the members of $\mathcal{A}$ are distinct. In addition, and for convenience only, we consider all samples from $\mathcal{A}$ as though drawn sequentially, that is, obtained with order.

Inequality (71) is imposed so that various expressions have simpler forms [see, e.g., (84)], in order to leave at least one unsampled individual with which to form an exchangeable pair, and also to yield

$$\lambda = \frac{N}{n(N-n)} \in (0, 1). \tag{74}$$

To form an exchangeable pair, let $X', X'', X_2, \ldots, X_n$ be a simple random sample of size $n + 1$ from $\mathcal{A}$, that is, with distribution

$$P(X' = x', X'' = x'', X_2 = x_2, \ldots, X_n = x_n)$$
$$= N_{n+1}^{-1} \mathbf{1}(\{x', x'', x_2, \ldots, x_n\} \subset \mathcal{A}, \text{ distinct}),$$

where $N_k = N!/(N - k)!$, the falling factorial. The pair

$$Y' = X' + \sum_{i=2}^n X_i \quad \text{and} \quad Y'' = X'' + \sum_{i=2}^n X_i$$

is clearly exchangeable with common marginal distribution that of $Y$ in (69). Since

$$E(X'|Y') = \frac{1}{n} Y' \quad \text{and} \quad E(X''|Y') = -\frac{1}{N - n} Y',$$

with $\lambda$ as in (74), we have

$$E(Y''|Y') = E(Y' - X' + X''|Y') = (1 - \lambda)Y',$$

proving that linearity condition (9) is satisfied.

We now follow the construction of the zero bias variable outlined in Section 3.3. Since $Y' - Y'' = X' - X''$, choose $X^\dagger, X^\ddagger$ independently of $X', X'', X_2, \ldots, X_n$, and with distribution proportional to the squared difference $(Y' - Y'')^2 = (X' - X'')^2$, that is, according to the distribution

$$(75) \qquad q(a, b) = \frac{(a - b)^2}{2N} \mathbf{1}(\{a, b\} \subset \mathcal{A}).$$

Now, the remainder of the sample from which we will construct $Y^\dagger$ and $Y^\ddagger$ must have the conditional distribution of $X_2, \ldots, X_n$ given $X' = X^\dagger, X'' = X^\ddagger$, that is, it must be a simple random sample of size $n - 1$ from $\mathcal{A} \setminus \{X^\dagger, X^\ddagger\}$.

However, we would like these $n - 1$ variables to correspond as closely as possible to the values in $\{X_2, \ldots, X_n\}$. For this reason, consider the difference and intersection

$$\mathscr{S} = \{X_2, \ldots, X_n\} \setminus \{X^\dagger, X^\ddagger\} \quad \text{and} \quad \mathcal{R}' = \{X_2, \ldots, X_n\} \cap \{X^\dagger, X^\ddagger\}.$$

The difference set $\mathscr{S}$ contains the variables in our original sample which can be used in the sample taken according to the conditional distribution given the inclusion of $X^\dagger$ and $X^\ddagger$, and $\mathcal{R}'$ contains the variables which cannot be common to both samples, that is, variables which must be replaced by others when forming $Y^\dagger$ and $Y^\ddagger$. In particular, if the intersection $\mathcal{R}'$ is empty, then $\{X_2, \ldots, X_n\}$ serves as the size $n - 1$ simple random sample from the complement of $\{X^\dagger, X^\ddagger\}$. Otherwise, $\mathcal{R}'$ is of size 1 or 2 and variables in $\mathcal{R}'$, in the order given by their indices, are replaced by those in a set $\mathcal{R}^\dagger$, of the same size as $\mathcal{R}'$, obtained by taking a simple random sample from the values available, that is, from the complement of

$$\mathcal{Q} = \{X_2, \ldots, X_n\} \cup \{X^\dagger, X^\ddagger\}.$$

In each case, the total resulting collection of the $n-1$ variables thus obtained are uniform from $\mathcal{A} \setminus \{X^\dagger, X^\ddagger\}$, that is, they have the conditional distribution of $X_2, \ldots, X_n$ given $X' = X^\dagger, X'' = X^\ddagger$; hence, (42) holds with

$$S = \sum_{a \in \mathscr{S}} a,$$

$$T' = \sum_{a \in \mathcal{R}'} a + X', \qquad T'' = \sum_{a \in \mathcal{R}'} a + X'',$$

$$T^\dagger = \sum_{a \in \mathcal{R}^\dagger} a + X^\dagger \quad \text{and} \quad T^\ddagger = \sum_{a \in \mathcal{R}^\dagger} a + X^\ddagger.$$

With $U \sim \mathcal{U}[0,1]$ independent of all other variables, by Proposition 3.2, a coupling of the zero biased variable $Y^*$ and $Y'$ is given by

$$Y^* = UX^\dagger + (1-U)X^\ddagger + S + \sum_{a \in \mathcal{R}^\dagger} a \quad \text{and} \quad Y' = X' + S + \sum_{a \in \mathcal{R}'} a,$$

and therefore their difference $V$ is given by

$$V = Y^* - Y' = UX^\dagger + (1-U)X^\ddagger - X' + \sum_{a \in \mathcal{R}^\dagger} a - \sum_{a \in \mathcal{R}'} a.$$

Now, using $X^\dagger =_d X^\ddagger$ and the independence of $U$, we may bound $E|V|$ by

$$(76) \qquad E|V| \leq E|X^\dagger| + E|X'| + E\left| \sum_{a \in \mathcal{R}^\dagger} a \right| + E\left| \sum_{a \in \mathcal{R}'} a \right|.$$

We bound the four terms of (76) separately.

Since $E(X')^2 = 1/N$, we have

$$E|\sqrt{N}X'| \leq E(\sqrt{N}X')^2 = 1 \leq (E|\sqrt{N}X'|^3)^{1/3} \leq E|\sqrt{N}X'|^3,$$

which gives the following bound on the second term of (76):

$$(77) \qquad \frac{1}{N}\sum_{a \in \mathcal{A}} |a| = E|X'| = \frac{1}{\sqrt{N}}E|\sqrt{N}X'| \leq \frac{1}{\sqrt{N}}E|\sqrt{N}X'|^3 = \sum_{a \in \mathcal{A}} |a|^3 = a_3.$$

From (75), the marginal distribution of $X^\dagger$ equals

$$q_1(a) = \frac{1}{2}\left(a^2 + \frac{1}{N}\right) \qquad \text{for } a \in \mathcal{A}.$$

Therefore, for the first term in (76), using (77), we have

$$(78) \qquad E|X^\dagger| = \sum_{a \in \mathcal{A}} |a| q_1(a) = \frac{1}{2}\sum_{a \in \mathcal{A}} |a|^3 + \frac{1}{2N}\sum_{a \in \mathcal{A}} |a| \leq a_3.$$

Moving to the last term in (76), since $\{X_2, \ldots, X_n\}$ and $X^\dagger, X^\ddagger$ are independent, for any $a \in \mathcal{A}$,

$$
\begin{aligned}
P(a \in \mathcal{R}') &= P(a \in \{X_2, \ldots, X_n\} \cap \{X^\dagger, X^\ddagger\}) \\
&= P(a \in \{X_2, \ldots, X_n\}) P(a \in \{X^\dagger, X^\ddagger\}) \\
&= 2P(a \in \{X_2, \ldots, X_n\}) P(X^\dagger = a) = \left(\frac{n-1}{N}\right)\left(a^2 + \frac{1}{N}\right),
\end{aligned}
$$

which implies that

$$
\begin{aligned}
E\left|\sum_{a \in \mathcal{R}'} a\right| &\leq E\sum_{a \in \mathcal{R}'} |a| = \sum_{a \in \mathcal{A}} |a| P(a \in \mathcal{R}') = \frac{n-1}{N}\sum_{a \in \mathcal{A}} |a|\left(a^2 + \frac{1}{N}\right) \\
&= \frac{n-1}{N}\left(\sum_{a \in \mathcal{A}} |a|^3 + \frac{1}{N}\sum_{a \in \mathcal{A}} |a|\right) \leq \frac{2n}{N}a_3,
\end{aligned}
\tag{79}
$$

using (77).

Beginning in a similar way for the third term in (76), since $P(|\mathcal{R}^\dagger| \in \{0, 1, 2\}) = 1$ and $P(a \in \mathcal{R}^\dagger, |\mathcal{R}^\dagger| = 0) = 0$ for all $a$, we have

$$
\begin{aligned}
E\left|\sum_{a \in \mathcal{R}^\dagger} a\right| &\leq \sum_{a \in \mathcal{R}^\dagger} |a| P(a \in \mathcal{R}^\dagger) \\
&= \sum_{a \in \mathcal{R}^\dagger} |a| P(a \in \mathcal{R}^\dagger, |\mathcal{R}^\dagger| = 1) + \sum_{a \in \mathcal{R}^\dagger} |a| P(a \in \mathcal{R}^\dagger, |\mathcal{R}^\dagger| = 2).
\end{aligned}
\tag{80}
$$

By independence, the joint distribution of $(X_2, \ldots, X_n)$ and $X^\dagger, X^\ddagger$, whose realizations are denoted $\chi_{n-1}$ and $u, v$, respectively, is given by

$$
p(\chi_{n-1}, u, v) = (N)_{n-1}^{-1}\mathbf{1}(\{x_2, \ldots, x_n\} \subset \mathcal{A}, \text{ distinct})q(u, v),
\tag{81}
$$

with $q(u, v)$ as in (75). Without further mention we consider only the event of probability one where $\chi_{n-1}$ is composed of distinct elements and $u \neq v$. Although $\chi_{n-1}$ is ordered, with a slight abuse of notation, we treat $\chi_{n-1}$ as an unordered set in expressions containing set operations, such as $\chi_{n-1} \cap \{u, v\}$. Taking $\mathcal{B}$ to be an ordered subset of $\mathcal{A}$ of size 1 or 2, the conditional distribution that $\mathcal{R}^\dagger = \mathcal{B}$, given $\chi_{n-1}$ and $u, v$, is uniform over all sets the size of the intersection of $\chi_{n-1}$ and $u, v$, taken from the complement of their union, that is,

$$
\begin{aligned}
p(\mathcal{B}|\chi_{n-1}, u, v) &= \frac{1}{(N - |\chi_{n-1} \cup \{u, v\}|)_{|\mathcal{B}|}} \\
&\quad \times \mathbf{1}(\mathcal{B} \cap (\chi_{n-1} \cup \{u, v\}) = \varnothing, |\mathcal{B}| = |\chi_{n-1} \cap \{u, v\}|).
\end{aligned}
$$

In particular, then, for $\mathcal{B}$ of size 1, using (81), we have

$$
\begin{aligned}
P(a \in \mathcal{R}^\dagger, |\mathcal{R}^\dagger| = 1) &= \sum_{\chi_{n-1}, u, v} p(a|\chi_{n-1}, u, v) p(\chi_{n-1}, u, v) \\
&= 2 \sum_{u \in \chi_{n-1}, v \notin \chi_{n-1}} p(a|\chi_{n-1}, u, v) p(\chi_{n-1}, u, v) \\
&= 2 \sum_{\substack{u \in \chi_{n-1}, v \notin \chi_{n-1} \\ \{u,v\} \cup \chi_{n-1} \not\ni a}} \frac{1}{N-n} \frac{1}{(N)_{n-1}} q(u,v) \\
&= \left( \frac{2}{N-n} \right) \frac{1}{(N)_{n-1}} \sum_{\{u,v\} \not\ni a} q(u,v) \sum_{\chi_{n-1} \cap \{v,a\} = \varnothing, \chi_{n-1} \ni u} 1
\end{aligned}
$$

$$
\tag{82} = \frac{2(n-1)(N-3)_{n-2}}{(N-n)(N)_{n-1}} \sum_{\{u,v\} \not\ni a} q(u,v)
$$

$$
= \frac{2(n-1)(N-n+1)}{(N)_3} \sum_{\{u,v\} \not\ni a} q(u,v)
$$

$$
\tag{83} = \frac{2(n-1)(N-n+1)}{(N)_3} \left( \left( 1 - \frac{1}{N} \right) - a^2 \right),
$$

where, in (82), the factor $(N-3)_{n-2}$ counts the number of ways that the $n-2$ additional elements required in $\chi_{n-1}$ can be taken from the $N-3$ available and the $n-1$ counts the number of positions that $u$ could occupy in the ordered set $\chi_{n-1}$. In addition, in the last equality, we have used

$$
\begin{aligned}
\sum_{\{u,v\} \not\ni a} q(u,v) &= \frac{1}{2N} \sum_{\{u,v\} \not\ni a} (u-v)^2 = \frac{1}{2N} \sum_{\{u,v\} \not\ni a} (u^2 - 2uv + v^2) \\
&= \frac{1}{N} \sum_{\{u,v\} \not\ni a} (u^2 - uv) = \frac{1}{N} \sum_{\{u,v\} \not\ni a} u^2 - \frac{1}{N} \sum_{\{u,v\} \not\ni a} uv \\
&= \frac{N-1}{N} \sum_{u \neq a} u^2 + \frac{a}{N} \sum_{u \neq a} u \\
&= \frac{N-1}{N} (1 - a^2) - \frac{a^2}{N} \\
&= \left( 1 - \frac{1}{N} \right) - a^2.
\end{aligned}
$$

Dropping the $-a^2$ term in (83) to get an upper bound and using (77) and the fact that $N \geq n \geq 3$, we have the following upper bound on the first

term in (80):

$$
(84) \qquad \sum_{a \in \mathcal{A}} |a| P(a \in \mathcal{R}^{\dagger}, |\mathcal{R}^{\dagger}| = 1) \leq \frac{2(n-1)(N-n+1)}{(N)_3} \left(1 - \frac{1}{N}\right) \sum_{a \in \mathcal{A}} |a|
$$

$$
\leq \frac{2(n-1)(N-n+1)}{N(N-2)} a_3 \leq \frac{2n}{N} a_3.
$$

To handle the second term in (80), we have, likewise, for $a$ and $b$ distinct,

$$
P(\mathcal{R}^{\dagger} = (a,b)) = \sum_{\chi_{n-1}, u, v} p(a, b | \chi_{n-1}, u, v) p(\chi_{n-1}, u, v)
$$

$$
= \sum_{\substack{|\{u,v\} \cap \chi_{n-1}| = 2 \\ (\{u,v\} \cup \chi_{n-1}) \cap \{a,b\} = \varnothing}} \frac{1}{(N-n+1)_2} \frac{1}{(N)_{n-1}} q(u, v)
$$

$$
= \frac{1}{(N-n+1)_2} \frac{1}{(N)_{n-1}} \sum_{\{u,v\} \cap \{a,b\} = \varnothing} q(u, v) \sum_{\substack{\chi_{n-1} \supset \{u,v\}, \\ \chi_{n-1} \cap \{a,b\} = \varnothing}} 1
$$

$$
= \frac{(n-1)_2}{(N-n+1)_2} \frac{(N-4)_{n-3}}{(N)_{n-1}} \sum_{\{u,v\} \cap \{a,b\} = \varnothing} q(u, v)
$$

$$
= \frac{(n-1)_2}{(N)_2} \frac{1}{(N-2)_2} \sum_{\{u,v\} \cap \{a,b\} = \varnothing} q(u, v)
$$

$$
= \frac{(n-1)_2}{(N)_2} \frac{1}{(N-2)_2} \frac{1}{N} ((N-2)(1 - a^2 - b^2) - (a+b)^2).
$$

Using symmetry, summing over $b \neq a$ and multiplying by 2 (since $a$ can be chosen as the first or second variable in the set $\mathcal{R}^{\dagger}$ of size 2) yields

$$
P(a \in \mathcal{R}^{\dagger}, |\mathcal{R}^{\dagger}| = 2) = \frac{2(n-1)_2}{(N)_2} \frac{1}{(N-2)_2} \frac{1}{N} ((N-1)(N-3) - (N^2 - 3N)a^2).
$$

By (71), $N > 3$, over which range the factor $-(N^2 - 3N)$ multiplying $a^2$ is negative; discarding it yields the upper bound

$$
P(a \in \mathcal{R}^{\dagger}, |\mathcal{R}^{\dagger}| = 2) \leq \frac{2(n-1)_2}{(N)_2} \frac{(N-1)(N-3)}{(N-2)_2} \frac{1}{N} = \frac{2(n-1)_2}{N^2(N-2)},
$$

so, by (77),

$$
(85) \qquad \sum_{a \in \mathcal{A}} |a| P(a \in \mathcal{R}^{\dagger}, |\mathcal{R}^{\dagger}| = 2) \leq \frac{2(n-1)_2}{N(N-2)} a_3 \leq 2 \left(\frac{n}{N}\right)^2 a_3.
$$

Inequalities (80), (84) and (85) yield the upper bound on the third term in (76),

$$(86) \qquad E\left|\sum_{a\in\mathcal{R}^\dagger} a\right| \le 2\left(\frac{n}{N} + \left(\frac{n}{N}\right)^2\right)a_3.$$

Combining the bounds on the four terms of (76) given in (77), (78), (79) and (86) gives

$$E|V| \le 2\left(1 + \frac{n}{N}\right)^2 a_3.$$

By (72), $EY = 0$, so $W = Y/\sigma$ and since $W^* = (Y/\sigma)^* = Y^*/\sigma$, Theorem 1.1 gives

$$\|F - \Phi\|_1 \le 2E|W^* - W| = \frac{2E|V|}{\sigma} \le \frac{4a_3}{\sigma}\left(1 + \frac{n}{N}\right)^2,$$

which is (73).   $\square$

**6. Combinatorial central limit theorem.**   We now use Theorem 1.1 to derive $L^1$ bounds for random variables $Y$ of the form

$$(87) \qquad Y = \sum_{i=1}^n a_{i,\pi(i)},$$

where $\pi$ is a permutation distributed uniformly over the symmetric group $\mathcal{S}_n$ and $\{a_{ij}\}_{1\le i,j\le n}$ are the components of a matrix $A \in \mathbf{R}^{n\times n}$. Letting

$$a_{\cdot\cdot} = \frac{1}{n^2}\sum_{i,j=1}^n a_{ij}, \qquad a_{i\cdot} = \frac{1}{n}\sum_{j=1}^n a_{ij} \quad\text{and}\quad a_{\cdot j} = \frac{1}{n}\sum_{i=1}^n a_{ij},$$

straightforward calculations show that the mean $\mu$ and variance $\sigma^2$ of $Y$ are given by

$$(88) \qquad \mu = na_{\cdot\cdot} \quad\text{and}\quad \sigma^2 = \frac{1}{n-1}\sum_{i,j}(a_{ij}^2 - a_{i\cdot}^2 - a_{\cdot j}^2 + a_{\cdot\cdot}^2);$$

the fact that (94) below is a probability distribution yields an equivalent representation for $\sigma^2$,

$$(89) \qquad \sigma^2 = \frac{1}{4n^2(n-1)}\sum_{i,j,k,l}[(a_{ik} + a_{jl}) - (a_{il} + a_{jk})]^2.$$

In what follows, we assume for the sake of nontriviality that $\sigma^2 > 0$. By (89), $\sigma^2 = 0$ if and only if $a_{il} - a_{i\cdot}$ does not depend on $i$, that is, if and only if the difference between any two rows of $A$ is some constant row vector.

Motivated by deriving approximating null distributions for permutation test statistics, Wald and Wolfowitz [25] proved the central limit theorem as $n \to \infty$ for the case where the factorization $a_{ij} = b_i c_j$ holds. This was later generalized by Hoeffding [12] to arrays $\{a_{ij}\}_{1 \leq i,j \leq n}$ in general. Motoo [17] gave Lindeberg-type sufficient conditions for the normal limit to hold.

In the supremum norm, von Bahr [2] and Ho and Chen [14] obtained Berry–Esseen bounds when the matrix $A$ is random, which yield the correct rate $O(n^{-1/2})$ only under some boundedness conditions. Bolthausen [6] obtained a bound of the correct order in terms of third-moment-type quantities, but with an unspecified constant. Goldstein [9] gave bounds of the correct order under boundedness, but with an explicit constant, for the cases where the random permutation $\pi$ is uniformly distributed and also when its distribution is constant on cycle type.

For each $n$, Theorem 6.1 provides an $L^1$ bound between the standardized variable $Y$ given in (87) and the normal, with an explicit constant depending on the third-moment-type quantity

$$(90) \qquad a_3 = \sum_{i,j=1}^n |a_{ij} - a_{i\cdot} - a_{\cdot j} + a_{\cdot\cdot}|^3.$$

When the elements of $A$ are all of comparable order, $\sigma^2$ is of order $n$ and $a_3$ of order $n^2$, making the bound below of order $n^{-1/2}$.

THEOREM 6.1.  *For $n \geq 3$, let $\{a_{ij}\}_{i,j=1}^n$ be the components of a matrix $A \in \mathbf{R}^{n \times n}$, let $\pi$ be a random permutation uniformly distributed over $\mathcal{S}_n$ and let $Y$ be given by* (87). *Then, with $\mu$, $\sigma^2$ and $a_3$ given in* (88) *and* (90), *$F$ the distribution function of $W = (Y - \mu)/\sigma$ and $\Phi$ that of the standard normal,*

$$\|F - \Phi\|_1 \leq \frac{a_3}{(n-1)\sigma^3} \left( 16 + \frac{56}{(n-1)} + \frac{8}{(n-1)^2} \right).$$

PROOF.  Since

$$Y - \mu = \sum_{i=1}^n (a_{i,\pi(i)} - a_{i\cdot} - a_{\cdot\pi(i)} + a_{\cdot\cdot}),$$

without loss of generality, we may replace $a_{ij}$ by $a_{ij} - a_{i\cdot} - a_{\cdot j} + a_{\cdot\cdot}$ in which case

$$(91) \qquad \sum_{i=1}^n a_{ij} = \sum_{j=1}^n a_{ij} = 0$$

and (90) becomes $a_3 = \sum_{ij} |a_{ij}|^3$. We will write $Y$ and $\pi$ interchangeably for $Y'$ and $\pi'$.

*Construction of $Y^\dagger, Y^\ddagger$:* We follow the construction outlined in Section 3.3; see also [9]. For $1 \le i, j \le n$, let $\tau_{ij}$ be the permutation which transposes $i$ and $j$. Given $\pi'$, take $(I, J)$ independent of $\pi'$, uniformly over all pairs $1 \le I \ne J \le n$, that is, with distribution

$$(92) \qquad p_1(i, j) = \frac{1}{(n)_2} \mathbf{1}(i \ne j).$$

Now, set $\pi'' = \pi' \tau_{I,J}$ and let $Y''$ be given by (87) with $\pi''$ replacing $\pi$. In particular, $\pi''(i) = \pi'(i)$ for $i \notin \{I, J\}$, so

$$(93) \qquad Y' - Y'' = (a_{I,\pi'(I)} + a_{J,\pi'(J)}) - (a_{I,\pi'(J)} + a_{J,\pi'(I)}).$$

We note that the difference depends only on $I, J, \pi'(I), \pi'(J)$ having distribution $p_1(i, j) p_1(k, l)$, where $k$ and $l$ are the realizations of $\pi'(I)$ and $\pi'(J)$, respectively. It can easily be shown (see [9]) that the pair $Y', Y''$ is exchangeable and satisfies the linearity condition (9) with $\lambda = 2/(n-1)$.

To construct $(Y^\dagger, Y^\ddagger)$ with distribution $(y' - y'')^2 dP(y', y'')/E(Y' - Y'')^2$ of (36), note first, using (93) and then (35) for the second equality, that

$$E(Y' - Y'')^2 = \frac{1}{n^2(n-1)^2} \sum_{i,j,k,l} [(a_{ik} + a_{jl}) - (a_{il} + a_{jk})]^2 = 2\lambda\sigma^2 = \frac{4\sigma^2}{n-1},$$

noting that the summand is zero if $i = j$ or $k = l$. Still following the outline given in Section 3.3, to begin the construction of $Y^\dagger$ and $Y^\ddagger$, choose $I^\dagger$, $J^\dagger$, $K^\dagger$, $L^\dagger$ independently of the remaining variables, according to their original difference distribution biased by the difference (93) squared, that is, with distribution

$$(94) \qquad \begin{aligned} p_2(i, j, k, l) &= \frac{[(a_{ik} + a_{jl}) - (a_{il} + a_{jk})]^2}{E(Y' - Y'')^2} p_1(i, j) p_1(k, l) \\ &= \frac{[(a_{ik} + a_{jl}) - (a_{il} + a_{jk})]^2}{4n^2(n-1)\sigma^2}; \end{aligned}$$

in particular, $P(I^\dagger = J^\dagger) = P(K^\dagger = L^\dagger) = 0$. Now, set

$$\pi^\dagger = \begin{cases} \pi \tau_{\pi^{-1}(K^\dagger), J^\dagger}, & \text{if } L^\dagger = \pi(I^\dagger), K^\dagger \ne \pi(J^\dagger), \\ \pi \tau_{\pi^{-1}(L^\dagger), I^\dagger}, & \text{if } L^\dagger \ne \pi(I^\dagger), K^\dagger = \pi(J^\dagger), \\ \pi \tau_{\pi^{-1}(K^\dagger), I^\dagger} \tau_{\pi^{-1}(L^\dagger), J^\dagger}, & \text{otherwise,} \end{cases}$$

and $\pi^\ddagger = \pi^\dagger \tau_{I^\dagger, J^\dagger}$. Note that $\{\pi^\dagger(I^\dagger), \pi^\dagger(J^\dagger)\} = \{\pi^\ddagger(I^\dagger), \pi^\ddagger(J^\dagger)\} = \{K^\dagger, L^\dagger\}$. As the conditional distribution of $\pi$, given that it takes particular values on some collection of indices, is uniform over all permutations restricted to take those values, the variables $Y^\dagger$ and $Y^\ddagger$ given by (87) with $\pi$ replaced by $\pi^\dagger$ and $\pi^\ddagger$, respectively, have joint distribution (36).

*Calculation of $E|Y^* - Y'|$*: By Proposition 3.2,

$$Y^* - Y' = UY^\dagger + (1-U)Y^\ddagger - Y'$$

$$= U \sum_{i=1}^n a_{i,\pi^\dagger(i)} + (1-U) \sum_{i=1}^n a_{i,\pi^\ddagger(i)} - \sum_{i=1}^n a_{i,\pi(i)}.$$

With

$$\mathcal{I} = \{I^\dagger, J^\dagger\} \cup \{\pi^{-1}(K^\dagger), \pi^{-1}(L^\dagger)\},$$

we see that if $i, j \notin \mathcal{I}$, then $\pi(i) = \pi^\dagger(j) = \pi^\ddagger(j)$.

Hence, setting $V = Y^* - Y'$, we have

$$(95) \qquad V = \sum_{i \in \mathcal{I}} (U a_{i,\pi^\dagger(i)} + (1-U) a_{i,\pi^\ddagger(i)} - a_{i,\pi(i)}).$$

Further, letting

$$R = |\{\pi(I^\dagger), \pi(J^\dagger)\} \cap \{K^\dagger, L^\dagger\}|$$

and $\mathbf{1}_k = \mathbf{1}(R = k)$, since $P(R \le 2) = 1$, we have

$$V = V\mathbf{1}_2 + V\mathbf{1}_1 + V\mathbf{1}_0 \quad \text{and therefore}$$

$$(96) \qquad E|V| \le E|V|\mathbf{1}_2 + E|V|\mathbf{1}_1 + E|V|\mathbf{1}_0.$$

The three terms on the right-hand side of (96) give rise to the three components of the bound in the theorem.

For notational simplicity, the following summations in this section are performed over all indices which appear, whether in the summands or in a (possibly empty) collection of restrictions. In what follows, we will have equalities and bounds such as

$$(97) \qquad \sum |a_{il}|[(a_{ik} + a_{jl}) - (a_{il} + a_{jk})]^2$$

$$= \sum |a_{il}|(a_{ik}^2 + a_{jl}^2 + a_{il}^2 + a_{jk}^2) \le 4n^2 a_3.$$

Due to the form of the square on the left-hand side, if the factors in a cross term agree in their first index, they will have differing second indices and likewise if their second indices agree. This gives cross terms which are zero by virtue of (91), since they will have at least one unpaired index outside the absolute value over which to sum, for instance, the index $k$ in the term $\sum |a_{il}| a_{ik} a_{il}$. Hence the equality. The inequality follows from the fact that for any choices $\iota_1, \iota_2, \kappa_1, \kappa_2 \in \{i, j, k, l\}$ with $\iota_1 \ne \kappa_1$ and $\iota_2 \ne \kappa_2$, perhaps by relabeling the indices appearing after the inequality,

$$(98) \qquad \sum_{i,j,k,l} |a_{\iota_1,\kappa_1}| a_{\iota_2,\kappa_2}^2 \le \left(\sum_{k,l} |a_{ij}|^3\right)^{1/3} \left(\sum_{i,j} |a_{kl}|^3\right)^{2/3} = n^2 a_3.$$

Generally, the power of $n$ in such an inequality, in this case 2, will be 2 less than the number of indices of summation, in this case 4.

*Decomposition on $R = 2$*: On $\mathbf{1}_2$, $\mathcal{I} = \{I^\dagger, J^\dagger\}$. As the intersection which gives $R = 2$ can occur in two different ways, we make the further decomposition

$$V\mathbf{1}_2 = V\mathbf{1}_{2,1} + V\mathbf{1}_{2,2},$$

where $\mathbf{1}_{2,1} = \mathbf{1}(\pi(I^\dagger) = K^\dagger, \pi(J^\dagger) = L^\dagger)$ and $\mathbf{1}_{2,2} = \mathbf{1}(\pi(I^\dagger) = L^\dagger, \pi(J^\dagger) = K^\dagger)$. Since $\pi^\dagger = \pi$ on $\mathbf{1}_{2,1}$, by (95),

$$
\begin{aligned}
V\mathbf{1}_{2,1} &= \sum_{i \in \{I^\dagger, J^\dagger\}} (U a_{i,\pi^\dagger(i)} + (1-U) a_{i,\pi^\ddagger(i)} - a_{i,\pi(i)}) \mathbf{1}_{2,1} \\
&= [U(a_{I^\dagger,\pi^\dagger(I^\dagger)} + a_{J^\dagger,\pi^\dagger(J^\dagger)}) \\
&\quad + (1-U)(a_{I^\dagger,\pi^\ddagger(I^\dagger)} + a_{J^\dagger,\pi^\ddagger(J^\dagger)}) - (a_{I^\dagger,\pi(I^\dagger)} + a_{J^\dagger,\pi(J^\dagger)})] \mathbf{1}_{2,1} \\
&= [U(a_{I^\dagger,\pi(I^\dagger)} + a_{J^\dagger,\pi(J^\dagger)}) \\
&\quad + (1-U)(a_{I^\dagger,\pi(J^\dagger)} + a_{J^\dagger,\pi(I^\dagger)}) - (a_{I^\dagger,\pi(I^\dagger)} + a_{J^\dagger,\pi(J^\dagger)})] \mathbf{1}_{2,1} \\
&= (1-U)(a_{I^\dagger,\pi(J^\dagger)} + a_{J^\dagger,\pi(I^\dagger)} - a_{I^\dagger,\pi(I^\dagger)} - a_{J^\dagger,\pi(J^\dagger)}) \mathbf{1}_{2,1} \\
&= (1-U)(a_{I^\dagger,L^\dagger} + a_{J^\dagger,K^\dagger} - a_{I^\dagger,K^\dagger} - a_{J^\dagger,L^\dagger}) \mathbf{1}_{2,1}.
\end{aligned}
\tag{99}
$$

Due to the presence of the indicator $\mathbf{1}_{2,1}$, taking the expectation of (99) requires a joint distribution which includes the values taken on by $\pi$ at $I^\dagger$ and $J^\dagger$, say $s$ and $t$, respectively. Since $s$ and $t$ can be any two distinct values and are independent of $I^\dagger, J^\dagger, K^\dagger$ and $L^\dagger$, we have, with $p_1$ and $p_2$ given in (92) and (94), respectively,

$$
\begin{aligned}
&p_3(i,j,k,l,s,t) \\
&\quad = P((I^\dagger, J^\dagger, K^\dagger, L^\dagger, \pi(I^\dagger), \pi(J^\dagger)) = (i,j,k,l,s,t)) \\
&\quad = p_2(i,j,k,l) p_1(s,t) = \frac{[(a_{ik} + a_{jl}) - (a_{il} + a_{jk})]^2}{4n^3(n-1)^2\sigma^2} \mathbf{1}(s \neq t).
\end{aligned}
\tag{100}
$$

Now, bounding the absolute value of the first term in (99) using (97),

$$
\begin{aligned}
E|(1-U) a_{I^\dagger,L^\dagger}| \mathbf{1}_{2,1} &= \frac{1}{2} \sum |a_{il}| \mathbf{1}(s = k, t = l) p_3(i,j,k,l,s,t) \\
&= \frac{1}{2} \sum |a_{il}| p_3(i,j,k,l,k,l) \\
&= \frac{1}{8n^3(n-1)^2\sigma^2} \sum |a_{il}| [(a_{ik} + a_{jl}) - (a_{il} + a_{jk})]^2 \\
&\leq \frac{a_3}{2n(n-1)^2\sigma^2}.
\end{aligned}
$$

Using the triangle inequality in (99) and applying the same reasoning to the remaining three terms shows that $E|V|\mathbf{1}_{2,1} \le 2a_3/(n(n-1)^2\sigma^2)$; since, by symmetry, the term $V\mathbf{1}_{2,2}$ can be bounded in this same way, we obtain

$$(101) \qquad E|V|\mathbf{1}_2 \le \frac{4a_3}{n(n-1)^2\sigma^2} \le \frac{4a_3}{(n-1)^3\sigma^2}.$$

*Decomposition on $R = 1$:* As the event $R = 1$ can occur in four different ways, depending on which element of $\{\pi(I^\dagger), \pi(J^\dagger)\}$ equals an element of $\{K^\dagger, L^\dagger\}$, we decompose $\mathbf{1}_1$ to yield

$$(102) \qquad V\mathbf{1}_1 = V\mathbf{1}_{1,1} + V\mathbf{1}_{1,2} + V\mathbf{1}_{1,3} + V\mathbf{1}_{1,4},$$

where $\mathbf{1}_{1,1} = \mathbf{1}(\pi(I^\dagger) = K^\dagger$ and $\pi(J^\dagger) \ne L^\dagger)$, on which $\mathcal{I} = \{I^\dagger, J^\dagger, \pi^{-1}(L^\dagger)\}$, specifying the remaining three indicators in (102) similarly. Now, using (95), and the fact that on $\mathbf{1}_{1,1}$, we have $\pi^\dagger = \pi\tau_{\pi^{-1}(L),J^\dagger}$ and $\pi^\ddagger = \pi\tau_{\pi^{-1}(L),J^\dagger}\tau_{J^\dagger,I^\dagger}$, so that $\pi^\dagger(\pi^{-1}(L)) = \pi^\ddagger(\pi^{-1}(L)) = \pi(J)$ it follows that

$$
\begin{aligned}
V\mathbf{1}_{1,1} &= \sum_{i \in \{I^\dagger, J^\dagger, \pi^{-1}(L^\dagger)\}} (Ua_{i,\pi^\dagger(i)} + (1-U)a_{i,\pi^\ddagger(i)} - a_{i,\pi(i)})\mathbf{1}_{1,1} \\
&= [U(a_{I^\dagger,\pi^\dagger(I^\dagger)} + a_{J^\dagger,\pi^\dagger(J^\dagger)} + a_{\pi^{-1}(L^\dagger),\pi^\dagger(\pi^{-1}(L^\dagger))}) \\
&\quad + (1-U)(a_{I^\dagger,\pi^\ddagger(I^\dagger)} + a_{J^\dagger,\pi^\ddagger(J^\dagger)} + a_{\pi^{-1}(L^\dagger),\pi^\ddagger(\pi^{-1}(L^\dagger))}) \\
&\qquad - (a_{I^\dagger,\pi(I^\dagger)} + a_{J^\dagger,\pi(J^\dagger)} + a_{\pi^{-1}(L^\dagger),\pi(\pi^{-1}(L^\dagger))})]\mathbf{1}_{1,1} \\
&= [U(a_{I^\dagger,K^\dagger} + a_{J^\dagger,L^\dagger} + a_{\pi^{-1}(L^\dagger),\pi(J^\dagger)}) \\
&\quad + (1-U)(a_{I^\dagger,L^\dagger} + a_{J^\dagger,K^\dagger} + a_{\pi^{-1}(L^\dagger),\pi(J^\dagger)}) \\
&\qquad - (a_{I^\dagger,K^\dagger} + a_{J^\dagger,\pi(J^\dagger)} + a_{\pi^{-1}(L^\dagger),L^\dagger})]\mathbf{1}_{1,1} \\
&= [Ua_{J^\dagger,L^\dagger} + (1-U)(a_{I^\dagger,L^\dagger} + a_{J^\dagger,K^\dagger} - a_{I^\dagger,K^\dagger}) \\
&\quad - a_{J^\dagger,\pi(J^\dagger)} - a_{\pi^{-1}(L^\dagger),L^\dagger} + a_{\pi^{-1}(L^\dagger),\pi(J^\dagger)}]\mathbf{1}_{1,1}.
\end{aligned}
\tag{103}
$$

For the first term in (103), dropping the restriction $t \ne l$ and summing over $t$ to obtain the first inequality and then applying (97) with $|a_{il}|$ replaced by $|a_{jl}|$, we obtain

$$
\begin{aligned}
EU|a_{J^\dagger,L^\dagger}|\mathbf{1}_{1,1} &= \frac{1}{2}\sum |a_{jl}|\mathbf{1}(s = k, t \ne l)p_3(i,j,k,l,s,t) \\
&\le \frac{1}{8n^2(n-1)^2\sigma^2}\sum |a_{jl}|[(a_{ik} + a_{jl}) - (a_{il} + a_{jk})]^2 \\
&\le \frac{a_3}{2(n-1)^2\sigma^2}.
\end{aligned}
\tag{104}
$$

The second, third and fourth terms in (103) result in the bound (104), with $|a_{jl}|$ replaced by $|a_{il}|, |a_{jk}|$ and $|a_{ik}|$, respectively, and applying corresponding forms of (97) on each gives

$$(105) \quad E|Ua_{J^\dagger,L^\dagger} + (1-U)(a_{I^\dagger,L^\dagger} + a_{J^\dagger,K^\dagger} - a_{I^\dagger,K^\dagger})|\mathbf{1}_{1,1} \le \frac{2a_3}{(n-1)^2\sigma^2}.$$

For the fifth term in (103), involving $a_{J^\dagger,\pi(J^\dagger)}$ without a uniform variable factor, we obtain

$$
\begin{aligned}
E|a_{J^\dagger,\pi(J^\dagger)}|\mathbf{1}_{1,1} &= \sum |a_{jt}|\mathbf{1}(s=k, t\ne l)p_3(i,j,k,l,s,t) \\
(106) \qquad &\le \frac{1}{4n^3(n-1)^2\sigma^2} \sum |a_{jt}|[(a_{ik}+a_{jl})-(a_{il}+a_{jk})]^2 \\
&\le \frac{a_3}{(n-1)^2\sigma^2}.
\end{aligned}
$$

Note that for the final inequality, although the sum being bounded is not of the form (97), having the index $t$, the same reasoning applies and that, moreover, the five indices of summation require that $n^2$ be replaced by $n^3$ in (98).

To handle the sixth term in (103), involving $a_{\pi^{-1}(L^\dagger),L^\dagger}$, we need the joint distribution

$$
\begin{aligned}
p_4&(i,j,k,l,s,t,u) \\
&= P((I^\dagger,J^\dagger,K^\dagger,L^\dagger,\pi(I^\dagger),\pi(J^\dagger),\pi^{-1}(L^\dagger)) = (i,j,k,l,s,t,u)),
\end{aligned}
$$

accounting for the value $u$ taken on by $\pi^{-1}(L^\dagger)$. If $l$ equals $s$ or $t$, then $u$ is already fixed at $i$ or $j$, respectively; otherwise, $\pi^{-1}(L^\dagger)$ is free to take any of the remaining available $n-2$ values, with equal probability. Hence, with $p_3$ given by (100), we deduce that

$$
p_4(i,j,k,l,s,t,u) = \begin{cases} p_3(i,j,k,l,s,t), & \text{if } (l,u)\in\{(s,i),(t,j)\}, \\ p_3(i,j,k,l,s,t)\dfrac{1}{n-2}, & \text{if } l\notin\{s,t\} \text{ and } u\notin\{i,j\}, \\ 0, & \text{otherwise.} \end{cases}
$$

Note, for example, that on $\mathbf{1}_{1,1}$, where $\pi(I^\dagger)=K^\dagger$ and $\pi(J^\dagger)\ne L^\dagger$, the value $u$ of $\pi^{-1}(L^\dagger)$ is neither $I^\dagger$ nor $J^\dagger$, so the second case above is the relevant one and the vanishing of the first sum on the third line of the following display is to be expected.

Now, calculating using the density $p_4$, for the sixth term in (103), we have

$$
\begin{aligned}
E|a_{\pi^{-1}(L^\dagger),L^\dagger}|\mathbf{1}_{1,1} \\
= \sum |a_{ul}|\mathbf{1}(s=k, t\ne l)p_4(i,j,k,l,s,t,u)
\end{aligned}
$$

$$= \sum_{t \neq l} |a_{ul}| p_4(i,j,k,l,k,t,u)$$

$$= \sum |a_{ik}| p_3(i,j,k,k,k,t) + \frac{1}{n-2} \sum_{l \notin \{k,t\}, u \notin \{i,j\}} |a_{ul}| p_3(i,j,k,l,k,t)$$

$$= \frac{1}{n-2} \sum_{l \neq t, u \notin \{i,j\}} |a_{ul}| p_2(i,j,k,l) p_1(k,t)$$

(107)
$$= \frac{1}{(n)_3} \sum_{t \notin \{l,k\}, u \notin \{i,j\}} |a_{ul}| p_2(i,j,k,l)$$

$$= \frac{1}{(n)_2} \sum_{u \notin \{i,j\}} |a_{ul}| p_2(i,j,k,l)$$

$$\leq \frac{1}{4n^3(n-1)^2 \sigma^2} \sum |a_{ul}| [(a_{ik} + a_{jl}) - (a_{il} + a_{jk})]^2$$

(108)
$$\leq \frac{a_3}{(n-1)^2 \sigma^2},$$

where the final inequality is achieved using (97) in the same way as for (106).

The computation for the seventh term in (103) begins in the same way as that for the sixth, yielding (107) with $a_{ut}$ replacing $a_{ul}$, so that

(109)
$$E|a_{\pi^{-1}(L^\dagger), \pi(J^\dagger)}| \mathbf{1}_{1,1}$$

$$= \frac{1}{(n)_3} \sum_{t \notin \{l,k\}, u \notin \{i,j\}} |a_{ut}| p_2(i,j,k,l)$$

$$\leq \frac{1}{4(n)_3 n^2 (n-1) \sigma^2} \sum |a_{ut}| [(a_{ik} + a_{jl}) - (a_{il} + a_{jk})]^2$$

$$\leq \frac{n^2 a_3}{(n)_3 (n-1) \sigma^2}$$

$$\leq \frac{3a_3}{(n-1)^2 \sigma^2},$$

where we have applied reasoning as in (97) and replaced $n^2$ by $n^4$ in (98) due to the sum over six indices.

Returning to (103) and adding the contribution (105) from the first four terms together with (106), (108) and (109) from the fifth, sixth and seventh, respectively, we obtain $E|V|\mathbf{1}_{1,1} \leq 7a_3/((n-1)^2 \sigma^2)$. Since, by symmetry, all four terms on the right-hand side of (102) can be handled in the same way as the first, we obtain the following bound on the event $R = 1$:

(110)
$$E|V|\mathbf{1}_1 \leq \frac{28a_3}{(n-1)^2 \sigma^2}.$$

*Decomposition on $R = 0$:* We have

$$\mathbf{1}_0 = \mathbf{1}(\pi(I^\dagger) \notin \{K^\dagger, L^\dagger\}, \pi(J^\dagger) \notin \{K^\dagger, L^\dagger\}),$$

$$\mathcal{I} = \{I^\dagger, J^\dagger, \pi^{-1}(K^\dagger), \pi^{-1}(L^\dagger)\},$$

and, from (95),

$$V\mathbf{1}_0 = \sum_{i \in \{I^\dagger, J^\dagger, \pi^{-1}(K^\dagger), \pi^{-1}(L^\dagger)\}} (U a_{i,\pi^\dagger(i)} + (1-U) a_{i,\pi^\ddagger(i)} - a_{i,\pi(i)}) \mathbf{1}_0$$

$$\begin{aligned}
(111) \quad &= [U(a_{I^\dagger, K^\dagger} + a_{J^\dagger, L^\dagger}) + (1-U)(a_{I^\dagger, L^\dagger} + a_{J^\dagger, K^\dagger}) \\
&\quad + a_{\pi^{-1}(K^\dagger), \pi(I^\dagger)} + a_{\pi^{-1}(L^\dagger), \pi(J^\dagger)} \\
&\quad - (a_{I^\dagger, \pi(I^\dagger)} + a_{J^\dagger, \pi(J^\dagger)} + a_{\pi^{-1}(K^\dagger), K^\dagger} + a_{\pi^{-1}(L^\dagger), L^\dagger})] \mathbf{1}_0.
\end{aligned}$$

Since the first four terms in (111) have the same distribution, we bound their contribution to $E|V|\mathbf{1}_0$, using (97), by

$$\begin{aligned}
(112) \quad 4EU|a_{I^\dagger, K^\dagger}|\mathbf{1}_0 &\leq 4EU|a_{I^\dagger, K^\dagger}| \\
&= 2\sum |a_{ik}| p_2(i,j,k,l) \\
&= \frac{1}{2n^2(n-1)\sigma^2} \sum |a_{ik}| [(a_{ik} + a_{jl}) - (a_{il} + a_{jk})]^2 \\
&\leq \frac{2a_3}{(n-1)\sigma^2}.
\end{aligned}$$

The sum of the contributions from the fifth and sixth terms of (111) can be bounded as

$$\begin{aligned}
2E|&a_{\pi^{-1}(L^\dagger), \pi(J^\dagger)}|\mathbf{1}_0 \\
&= 2 \sum_{s \notin \{k,l\}, t \notin \{k,l\}} |a_{ut}| p_4(i,j,k,l,s,t,u) \\
&= \frac{2}{n-2} \sum_{s \notin \{k,l\}, t \notin \{k,l\}, u \notin \{i,j\}, s \neq t} |a_{ut}| p_3(i,j,k,l,s,t) \\
(113) \quad &\leq \frac{n-3}{2(n-2)n^3(n-1)^2\sigma^2} \sum |a_{ut}| [(a_{ik} + a_{jl}) - (a_{il} + a_{jk})]^2 \\
&\leq \frac{2n(n-3)a_3}{(n-2)(n-1)^2\sigma^2} \\
(114) \quad &\leq \frac{2a_3}{(n-1)\sigma^2},
\end{aligned}$$

where inequality (113) is obtained by summing over the $n-3$ choices of $s$ and dropping the remaining restrictions, and the next by following the reasoning of (97).

For the sum of the contributions from the seventh and eighth terms of (111), summing over the $n-3$ choices of $t$ and then dropping the remaining restrictions to obtain the first inequality, we have

$$
\begin{aligned}
2E|a_{I^{\dagger},\pi(I^{\dagger})}|\mathbf{1}_0 &= 2\sum_{s\notin\{k,l\},t\notin\{k,l\}}|a_{is}|p_3(i,j,k,l,s,t)\\
&= \frac{1}{2n^3(n-1)^2\sigma^2}\sum_{s\notin\{k,l\},t\notin\{k,l\},s\neq t}|a_{is}|[(a_{ik}+a_{jl})-(a_{il}+a_{jk})]^2\\
&\leq \frac{n-3}{2n^3(n-1)^2\sigma^2}\sum|a_{is}|[(a_{ik}+a_{jl})-(a_{il}+a_{jk})]^2\\
&\leq \frac{2(n-3)a_3}{(n-1)^2\sigma^2}\leq\frac{2a_3}{(n-1)\sigma^2}.
\end{aligned}
\tag{115}
$$

The total contribution of the ninth and tenth terms together can be bounded like the sum of the fifth and sixth, yielding (113) with $|a_{ul}|$ replacing $|a_{ut}|$, then summing over the $n$ choices of $t$ gives

$$
\begin{aligned}
&2E|a_{\pi^{-1}(L^{\dagger}),L^{\dagger}}|\mathbf{1}_0\\
&\leq \frac{n-3}{2(n-2)n^2(n-1)^2\sigma^2}\sum|a_{ul}|[(a_{ik}+a_{jl})-(a_{il}+a_{jk})]^2\\
&\leq \frac{2n(n-3)a_3}{(n-2)(n-1)^2\sigma^2}\leq\frac{2a_3}{(n-1)\sigma^2}.
\end{aligned}
\tag{116}
$$

Adding up the bounds for the first four terms (112), the fifth and sixth terms (114), the seventh and eighth terms (115) and the ninth and tenth terms (116) yields

$$
E|V|\mathbf{1}_0\leq\frac{8a_3}{(n-1)\sigma^2}.
\tag{117}
$$

Since $W^*=(Y/\sigma)^*=Y^*/\sigma$, we have $E|W^*-W|=E|V|/\sigma$. Hence, summing the $R=2$, $R=1$ and $R=0$ contributions to $E|V|$ given in (101), (110) and (117), respectively, the proof of the theorem is completed by applying Theorem 1.1. $\square$

## 7. Remarks.
In Section 3.2, a new method of constructing zero bias couplings is presented which closely parallels the construction for size bias couplings. Applying also an existing construction, the zero bias method for computing $L^1$ bounds to the normal is illustrated in four situations.

The zero bias transformation for normal approximation is not restricted to the $L^1$ norm. The supremum norm is considered in [9] through the use of smoothing inequalities, although useful bounds there are only obtained when $|Y^*-Y|$ can be almost surely bounded by a quantity small relative to $\mathrm{Var}(Y)$. This restriction at present prevents the application of the zero bias method from computing supremum norm bounds in various examples, cone measure being one. It is hoped that this restriction may be relaxed in future work.

The Stein equation also presents the possibility for deriving total variation bounds in a way similar to the manner in which the $L^1$ bounds used here were derived in [8]. Letting a random variable denote its own distribution, recall that the total variation distance between the distributions of $X$ and $Y$ can be defined in terms of differences in expectations over bounded measurable test functions $h$:

$$(118) \qquad \|X - Y\|_{\mathrm{TV}} = \tfrac{1}{2} \sup_{|h| \leq 1} |Eh(X) - Eh(Y)|.$$

Now, consider the Stein equation, with $\sigma^2 = 1$, say, for such an $h$,

$$f'(x) - xf(x) = h(x) - Eh(Z),$$

where $Z$ is a standard normal variable. Stein [23] shows that if $|h| \leq 1$, then $f$ is differentiable with $|f'| \leq 2$ and hence, for a mean zero variance 1 random variable $W$,

$$\begin{aligned}
|Eh(W) - Eh(Z)| &= |Ef'(W) - EWf(W)| \\
&= |Ef'(W) - Ef'(W^*)| \leq 4\|W - W^*\|_{\mathrm{TV}}.
\end{aligned}$$

Dividing by 2 and taking supremum over $h$ as indicated in (118) yields

$$\|W - Z\|_{\mathrm{TV}} \leq 2\|W - W^*\|_{\mathrm{TV}},$$

a total variation bound parallel to the $L^1$ bound in Theorem 1.1.

**Acknowledgments.** The author would like to sincerely thank a reviewer for many helpful suggestions, in particular, regarding the formulation of Proposition 2.2 and the computation of $c_1$, and for helping to simply the argument and reduce the magnitude of the constants in Theorem 6.1.

Department of Mathematics
University of Southern California
Kaprielian Hall, Room 108
3620 Vermont Avenue
Los Angeles, California 90089-2532
USA
E-mail: larry@math.usc.edu